\documentclass[preprint,11pt]{elsarticle}

\usepackage[a4paper]{geometry}
\usepackage{amssymb}
\usepackage{amsfonts}
\usepackage{graphicx}
\usepackage{caption}
\usepackage{pgf}
\usepackage{subfig}
\usepackage{mathtools} 
\usepackage{acro} 
\usepackage{comment} 
\usepackage[ruled, vlined, english, onelanguage, algo2e]{algorithm2e}
\usepackage{xfrac} 
\usepackage{makecell} 
\usepackage{hyperref}
\usepackage[capitalise]{cleveref}
\Crefname{equation}{}{}

\newcommand{\orcid}[1]{\href{https://orcid.org/#1}{\includegraphics[width=10pt]{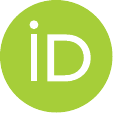}}}

\numberwithin{equation}{section}

\usepackage{fancyhdr}
\newcommand*\eps{\varepsilon} 
\newcommand*\nid{\text{d} }   
\newcommand{\pd}[2]{\frac{\partial #1}{\partial #2} } 
\newcolumntype{?}[1]{!{\vrule width #1pt}} 

\newtheorem{thm}{Theorem}

\newdefinition{rmk}{Remark}
\newproof{pf}{Proof}

\Crefname{algocf}{Algorithm}{Algorithms}

\DeclareAcronym{ode}{
	short = ODE,
	long  = ordinary differential equation
}
\DeclareAcronym{odes}{
	short = ODEs,
	long  = ordinary differential equations
}

\DeclareAcronym{pde}{
	short = PDE,
	long  = partial differential equation
}
\DeclareAcronym{pdes}{
	short = PDEs,
	long  = partial differential equations
}

\DeclareAcronym{ssp}{
	short = SSP,
	long  = strong stability preserving
}

\DeclareAcronym{tvd}{
	short = TVD,
	long  = total varitation diminishing
}

\DeclareAcronym{dg}{
	short = DG,
	long  = discontinuous Galerkin
}

\DeclareAcronym{dgsem}{
	short = DGSEM,
	long  = discontinuous Galerkin spectral element method
}

\DeclareAcronym{hllc}{
	short = HLLC,
	long  = Harten-Lax-Van Leer Contact
}

\DeclareAcronym{hll}{
	short = HLL,
	long  = Harten-Lax-Van Leer
}

\DeclareAcronym{mhd}{
	short = MHD,
	long  = magnetohydrodynamics
}

\DeclareAcronym{Mhd}{
	short = MHD,
	long  = Magnetohydrodynamics
}

\DeclareAcronym{glm-mhd}{
	short = GLM-MHD,
	long  = generalized Lagrangian multiplier magnetohydrodynamics
}

\DeclareAcronym{PERK}{
	short = P-ERK,
	long = Paired-Explicit Runge-Kutta
}

\journal{Springer Journal of Scientific Computing}

\begin{document}

	\begin{frontmatter}
		
		\title{Many-Stage Optimal Stabilized Runge-Kutta Methods for Hyperbolic Partial Differential Equations}
		
		\cortext[cor1]{Corresponding author}
		
		\author[1]{Daniel Doehring\corref{cor1} \orcid{0009-0005-4260-0332}}
		\author[2]{Gregor J. Gassner \orcid{0000-0002-1752-1158}}
		\author[1]{and Manuel Torrilhon \orcid{0000-0003-0008-2061}}
		
		\affiliation[1]{
			organization={Applied~and~Computational~Mathematics,~RWTH~Aachen~University},
			addressline={\\Schinkelstrasse~2}, 
			city={Aachen},
			postcode={52062}, 
			state={North~Rhine-Westphalia},
			country={Germany}.
		}
		\affiliation[2]{
			organization={Department of Mathematics and Computer Science, Center for Data and Simulation Science, \\University of Cologne},
			addressline={Weyertal 86-90},
			city={Cologne},
			postcode={50931},
			state={North~Rhine-Westphalia},
			country={Germany}.
		}
		
		\begin{abstract}
			A novel optimization procedure for the generation of stability polynomials of stabilized explicit Runge-Kutta methods is devised.
			Intended for semidiscretizations of hyperbolic partial differential equations, the herein developed approach 
			allows the optimization of stability polynomials with more than hundred stages.
			A potential application of these high degree stability polynomials are problems with locally varying characteristic speeds as found in non-uniformly refined meshes and different wave speeds.
			
			To demonstrate the applicability of the stability polynomials we construct $2N$ storage many-stage Runge-Kutta methods that match their designed second order of accuracy when applied to a range of linear and nonlinear hyperbolic PDEs with smooth solutions.
			The methods are constructed to reduce the amplification of round off errors which becomes a significant concern for these many-stage methods.
		\end{abstract}
		
		
	\end{frontmatter}
	
	
	%
	\section{Stabilized Explicit Runge-Kutta Methods}
	Explicit Runge-Kutta methods are commonly considered the default choice for the integration of hyperbolic \ac{pdes}.
	In contrast to implicit methods, explicit methods require only evaluations of the right-hand side instead of solving a potentially nonlinear system of equations.
	In fact, in the context of hyperbolic \ac{pdes} typically nonlinear fluxes are of interest.
	Explicit methods come with the drawback that the maximum stable timestep needs to be significantly reduced which is in the context of (hyperbolic) \ac{pdes} commonly referred to as the CFL condition \cite{courant1967partial}.
	
	To increase computational efficiency, \textit{stabilized} explicit Runge-Kutta methods have been introduced already in the 1960s \cite{franklin1958numerical, guillou1960domaine, saul1960integration} targeting semidiscretizations 
	\begin{subequations}
		\label{eq:ODESys}
		\begin{align}
			\boldsymbol U(t_0) &= \boldsymbol U_0 \\
			\label{eq:ODESys2}
			\boldsymbol U'(t) &= \boldsymbol F\big(\boldsymbol U(t) \big)
		\end{align}
	\end{subequations}
	of parabolic \ac{pdes}.
	The approach of discretizing the spatial derivatives and leaving time continuous at first is commonly referred to as the \textit{Method of Lines} approach.
	Although originally termed for finite difference based spatial discretizations \cite{forsythe1960finite} the name persisted and is used nowadays to comprise different discretization techniques like finite elements, finite volumes, and discontinuous Galerkin.
	
	The central idea of stabilized explicit Runge-Kutta methods is to use additional stages primarily to improve the stability properties of the scheme, i.e., allow for larger timesteps.
	In particular, one usually settles for a moderate order of accuracy and uses the additional degrees of freedom to improve the stability properties of the scheme for a certain spectrum \cite{wanner1996solving, hundsdorfer2003numerical, verwer1996explicit}.
	The fact that the eigenvalues of the Jacobian 
	\begin{equation}
		\label{eq:JacobianSemiDisc}
		J(\boldsymbol U) \coloneqq \pd{\boldsymbol{F}}{\boldsymbol{U}}
	\end{equation}
	are for parabolic \ac{pdes} distributed on the negative real axis enables a successful construction of stabilized explicit Runge-Kutta methods optimized for this special case, see for instance the reviews \cite{abdulle2015explicit, van1996development}.
	The fact that the optimal stability polynomial $P_{S,1}(z)$ of first order accurate methods with $S$ stages is given by the shifted Chebyshev polynomial of first kind $T_S(z)$ \cite{franklin1958numerical, guillou1960domaine, dinsome1958, burrage1985order}
	serves as a valuable guidance to construct higher order approximately optimal stability polynomials \cite{lomax1968construction, lebedev1993new1, lebedev1993new2, van1977construction, van1971numerical, abdulle2001second} which stand on solid ground due to proven existence and uniqueness of optimal stability polynomials with maximum (negative) real axis inclusion \cite{riha1972optimal}.
	
	For the purely hyperbolic case, i.e., for eigenvalues $\lambda \in \sigma(J)$ exclusively on the imaginary axis also results for first order accurate methods in terms of Chebyshev polynomials of first kind are known \cite{van1972explicit, van1977construction, kinnmark1984one1, sonneveld1984minimax}.
	Approximations to higher order accurate stability polynomials can be found in \cite{kinnmark1984one2}.
	
	For the more general case, i.e., where not only either the real or complex line are of interest, but rather a two-dimensional part of the complex plane, concrete results for stability polynomials are rare.
	To the best of our knowledge, only for the circle results for first \cite{jeltsch1978largest} and second \cite{owren1990some} order accurate optimal stability polynomials of variable degree $S$ are available.
	In contrast, one usually has to resort to numerical optimization of the stability polynomials.
	A nonextensive list includes optimized methods for certain geometrically defined spectra \cite{niegemann2012efficient, torrilhon2007essentially}, particular equations \cite{KENNEDY2000177, ALLAMPALLI20093837, TOULORGE20122067, MEAD1999404} and spatial discretization techniques \cite{al2021optimized, kubatko2014optimal, Nasab2022Optimal}.
	
	A general framework for maximizing the region of absolute stability of an explicit Runge-Kutta methods for a particular spectrum has been developed in \cite{ketcheson2013optimal} which has been extensively used \cite{VERMEIRE201955, kubatko2014optimal, SCHLOTTKELAKEMPER2021110467, VERMEIRE2019465, HEDAYATINASAB2022111470, VERMEIRE2021110022, hedayati2021optimal, klinge2018strong, al2021optimized, ellison2021parallel}.
	The approach presented in \cite{ketcheson2013optimal} is, however, for general spectra limited to polynomials of moderate $(16-20)$ degree due to issues with floating point precision.
	This is a consequence of the formulating the task to find the maximum admissible timestep as a convex optimization problem in terms of the monomial coefficients of the stability polynomial.
	The monomial coefficients are also the optimization variables in \cite{niegemann2012efficient, ALLAMPALLI20093837, TOULORGE20122067} where polynomials are optimized up to 14, 7, and 8 degrees, respectively.
	
	In this work, we develop a formulation that avoids the inherent ill-conditioning of the monomials-based approach at the cost of turning the optimization problem into a nonlinear one.
	In particular, instead of parametrizing the polynomial in terms of the monomial coefficients, the stability polynomial $P(z)$ is characterized by the roots of $\big(P(z)-1\big)/z$.
	This can be seen as a generalization of the approach employed in \cite{torrilhon2007essentially} wherein the stability polynomial is described by its extrema.
	While this seems daunting, it will be shown that for a relevant class of spectra an excellent initial guess can be supplied allowing the successful optimization of the highly nonlinear problem.
	By doing so, we are able to optimize stability polynomials of degrees larger than 100, which to the best of our knowledge is an unprecedented success for general spectra.
	
	A potential application of these high degree stability polynomials are the recently published \ac{PERK} methods \cite{VERMEIRE2019465, HEDAYATINASAB2022111470, VERMEIRE2023112159} which achieve local time stepping effects by combining a set of stabilized methods.
	In particular, due to their construction as partitioned Runge-Kutta methods \cite{HairerWanner1} they ensure consistency and conservation which are not trivially satisfied by other classical multirate methods \cite{Hundsdorfer2013}.
	Being based on the optimization approach developed in \cite{ketcheson2013optimal}, the stability polynomials and corresponding \ac{PERK} methods are currently limited to 16 stages.
	The availability of higher degree polynomials allows for a potentially even more efficient treatment of e.g. locally refined meshes or varying characteristic speeds.
	In this work, however, we focus on the optimization of high degree stability polynomials and methods that can be directly constructed thereof.
	The application to \ac{PERK} methods is left for future work.
	
	The paper is organized as follows. The motivation for the new optimization approach is presented in \cref{sec:preliminaries} and \cref{sec:motivation}, based on findings for the proven optimal stability polynomials for disk-like spectra.
	The necessary details required for successful optimization of the problem are discussed in \cref{sec:Formulation} which is followed by a discussion of the implementation in \cref{sec:Implementation}.
	In \cref{sec:OptStabPolysStrictlyConvexSpectra} a couple of high degree optimized stability polynomials are presented.
	Extensions of the original approach to non-convex spectra are covered in \cref{sec:NonConvexSpectra} and applied in \cref{sec:OptStabPolysNonConvex}.
	The construction of actual Runge-Kutta methods from the high-degree stability polynomials follows in \cref{sec:Construction} with a special focus on internal stability.
	\cref{sec:ResultsStrictlyConvexSpectra} presents the application of the many stage methods to linear and nonlinear problems.
	\cref{sec:Conclusions} concludes the paper.
	\section*{Note on Terminology}
	For the sake of readability, we will refer to the shifted Chebyshev polynomials of first kind simply as Chebyshev polynomials.
	As the extreme points of the (shifted) Chebyshev polynomials (of first kind) play a crucial role in this work, we will call them in the interest of brevity simply Chebyshev extreme points.
	Furthermore, we will refer to the the Chebyshev extreme points with extremal value $+1$ as positive Chebyshev extreme points.
	
	Since a significant part of the paper centers around stability polynomials from which Runge-Kutta methods of certain order may be constructed, we will say that a polynomial of degree $S$ is of order/accuracy $p$ whenever the first $i= 0, \dots, p$ coefficients match the first $i=0, \dots, p$ coefficients of the Taylor series of the exponential.
	\section{Preliminaries}
	\label{sec:preliminaries}
	Runge-Kutta methods are single step methods, i.e., when applied to the constant coefficient linear initial value problem (IVP)
	\begin{subequations}
		\label{eq:LinearConstantCoeffODEFuncIV}
		\begin{align}
			\boldsymbol u(t_0 = 0) &= \boldsymbol u_0 \\
			\label{eq:LinearConstantCoeffODEFunc}
			\boldsymbol u'(t) &= A \boldsymbol u(t)
		\end{align}
	\end{subequations}
	the new approximation $\boldsymbol u_{n+1}$ can be computed from the previous iterate $\boldsymbol u_{n}$ as
	\begin{equation}
		\label{eq:StabilityFunction}
		\boldsymbol u_{n+1} = R(\Delta t A) \boldsymbol u_{n}.
	\end{equation}
	Here, $R(z)$ denotes the \textit{stability function} \cite{wanner1996solving} of the Runge-Kutta method.
	For implicit methods $R(z)$ is a rational function, while for explicit methods $R(z)$ is a polynomial with real coefficients
	\begin{equation}
		\label{eq:PolynomialMonCoeffs}
		P_S(z; \boldsymbol \alpha ) = \sum_{j=0}^S \alpha_j z^j, \quad \boldsymbol \alpha \in \mathbb R^{S+1}.
	\end{equation}
	By comparing \cref{eq:StabilityFunction} to the solution $\boldsymbol u(t_{n+1}) = \exp(A \Delta t) \boldsymbol u(t_n)$ of \cref{eq:LinearConstantCoeffODEFunc} it follows from the definition of the exponential that for a $p'$th (linearly) consistent approximation the coefficients $\alpha_j$ need to satisfy 
	\begin{equation}
		\label{eq:ConsistencyRequirement}
		\alpha_j \overset{!}{=} \frac{1}{j!}, \quad j =0, \dots, p.	
	\end{equation}
	It is customary to define the family of polynomials with real coefficients over the complex numbers of degree $S$ and corresponding order of accuracy $p$ as $\mathcal{P}_{S,p}$.
	
	Due to the fact that the stability function $R(z)$ is a polynomial $P(z)$ for explicit methods it follows that the \textit{region of absolute stability} 
	\begin{equation}
		\mathcal S \coloneqq \{ z \in \mathbb C : \vert R(z) \vert \leq 1\}
	\end{equation}
	is necessarily bounded.
	We further define the \textit{boundary of the region of absolute stability} as 
	\begin{equation}
		\partial \mathcal S \coloneqq \{ z \in \mathbb C : \vert R(z) \vert = 1\}
	\end{equation}
	which will also be called in brief \textit{stability boundary}.
	\subsection{Optimization Objective}
	The optimization objective is now to find the maximum possible timestep $\Delta t_{S,p}^\star$ for a polynomial of degree $S$ corresponding to a method with order of accuracy $p$:
	\begin{equation}
		\label{eq:OptProb}
		\underset{P_{S,p} \in \mathcal{P}_{S,p}}{\max} \Delta t \text{ such that } \big \vert P_{S,p}(\Delta t \lambda^{(m)}) \big \vert \leq 1, \quad  m = 1 , \dots , M.
	\end{equation}
	Here, $\{ \lambda^{(m)} \}_{m=1, \dots, M} = \sigma(J)$ are the eigenvalues of the Jacobian \cref{eq:JacobianSemiDisc} where the \ac{ode} system \cref{eq:ODESys} corresponds in this work to the semidiscretization of \ac{pdes} describing typically physical processes.
	As a consequence, we assume that there are no amplifiying modes among the eigenvalues $\lambda^{(m)}$ corresponding to generation of energy, thus all eigenvalues should have non-positive real part:
	\begin{equation}
		\text{Re}\big(\lambda \big) \leq 0 \qquad \forall \: \lambda \in \sigma(J).
	\end{equation}
	It should be stressed that the maximum possible timestep $\Delta t^\star$ is the single optimization target in this work.
	In particular, we do not focus on the reduction of dispersion or dissipation errors as done for instance in \cite{hu1996low, berland2006low, bernardini2009general} or other objectives like maximum \ac{ssp} coefficients \cite{spiteri2002new, ruuth2006global, kubatko2014optimal}.
	\subsubsection{Convex Problem Formulation}
	As noted in \cite{ketcheson2013optimal} \cref{eq:OptProb} is for fixed timestep $\Delta t$ a convex optimization problem when parametrizing the stability polynomial $P_{S,p}$ in terms of the monomial coefficients $\boldsymbol \alpha$:
	\begin{equation}
		\label{eq:OptProbConvex}
		\underset{\boldsymbol \alpha \in \mathbb R^{S-p}}{\max} \Delta t \text{ such that } \big \vert P_{S,p}(\Delta t \lambda^{(m)}; \boldsymbol \alpha) \big \vert \leq 1, \quad  m = 1 , \dots , M.
	\end{equation}
	Paired with an outer bisection to determine the timestep, \cref{eq:OptProbConvex} leads to an efficient optimization routine that can be solved with standard software such as \texttt{SeDuMi} \cite{sturm1999using} or \texttt{ECOS} \cite{6669541} for second order cone programs.
	The downside of formulation \cref{eq:OptProbConvex} is that the monomial coefficients scale essentially as $\alpha \sim 1/(j!)$ which limits this approach in standard double precision to about 16 to 20 stages for general spectra.
	Relieving this inherent issue is the key 
	achievement
	of our approach motivated in the next subsection. 
	\subsection{Central Observation}
	The linear consistency requirement \cref{eq:ConsistencyRequirement} implies that every at least first order linearly consistent stability polynomial $P_{S,1}(z)$ can be written as
	\begin{equation}
		\label{eq:DefinitionLowerDegree}
		P_{S,1}(z; \widetilde{\boldsymbol r}) = 1 + z \prod_{j=1}^{S-1} \left( 1- \frac{z}{\widetilde{r}_j} \right) \\
		\eqqcolon 1 + z \widetilde{P}_{S-1}(z; \widetilde{\boldsymbol r})
	\end{equation}
	where the \textit{lower degree polynomial} $\widetilde{P}_{S-1}$ is parametrized in its complex-conjugated roots $\widetilde{\boldsymbol r} \in \mathbb C^{S-1}$.
	Clearly, the roots of the lower-degree polynomial (called from now on \textit{pseudo-extrema} of the original polynomial) form a subset of the stability boundary $\partial \mathcal S$:
	\begin{equation}
		\label{eq:PseudoExtremaSubsetStabBnd}
		\vert P_{S,1}(\widetilde{r}_j; \widetilde{\boldsymbol r})	\vert = \vert 1 + \widetilde{r}_j \underbrace{\widetilde{P}_{S-1}(\widetilde{r}_j; \widetilde{\boldsymbol r})}_{=0} \vert = 1 \quad \Rightarrow \big \{\widetilde{r}_j \big \}_{j=1, \dots, S-1} \subset \partial \mathcal S.
	\end{equation}
	This observation is particularly useful when the shape of the stability boundary $\partial \mathcal S$ can be a-priori roughly estimated.
	For optimized polynomials of higher degrees, $\partial \mathcal S$ will follow the spectrum closely, thus one can use an envelope of the spectrum itself as a reasonable approximation to $\partial \mathcal S$.
	This is discussed in \cref{sec:CentralAssumption} and \cref{sec:NonConvexSpectra} in more detail.

	\section{Motivating the Optimization in Pseudo-Extrema}
	\label{sec:motivation}
	The optimization problem \cref{eq:OptProb} for stability polynomials in terms of complex-conjugated pseudo-extrema reads
	\begin{equation}
		\label{eq:OptProbPE}
		\underset{\widetilde{\boldsymbol r} \in \mathbb C^{S-1}}{\max} \Delta t \text{ such that } \big \vert P_{S,p}(\Delta t \lambda^{(m)}; \widetilde{\boldsymbol r}) \big \vert \leq 1, \quad  m = 1 , \dots , M.
	\end{equation}
	This is a highly nonlinear optimization problem for which we need an excellent initial guess in order to have a reasonable chance to reach the global optimum.
	In the next section, we will thoroughly motivate a suitable assumption on the distribution of the pseudo-extrema based on results of proven optimal stability polynomials.
	\subsection{Pseudo-Extrema of Proven Optimal Stability Polynomials for Disks}
	\label{sec:PEProvenOptimalStabPolys}
	We direct our attention to examining the pseudo-extrema of known optimal stability polynomials $P_{S,p}(z)$ for circular spectra.
	For the sake of completeness, the case for parabolic spectra, i.e., eigenvalues on the negative real line, is provided in \ref{sec:PECheby}.
	\subsubsection{First Order Consistent Stability Polynomial}
	As shown in \cite{jeltsch1978largest}, the optimal first order stability polynomial with largest disk inclusion is given by a sequence of Forward Euler steps:
	\begin{equation}
		\label{eq:StabilityPolynomialDisk}
		P^\text{Disk}_{S, 1}(z) = \bigg( 1 + \frac{z}{S} \bigg)^S.
	\end{equation}
	For this stability polynomial it holds that the disk with radius $S$ 
	\begin{equation}
		D_S \coloneqq \{ z \in \mathbb C : \vert z + S \vert \leq S\}
	\end{equation}
	is contained in the region of absolute stability $\mathcal S$.
	Consider the positive Chebyshev extreme points $x_j \in \mathbb R: T_S(x_j) = 1$
	of the $S$ degree polynomial on the here relevant $[-2S, 0]$ interval which are given by
	\begin{equation}
		\label{eq:ChebyExtremrePointsReal}
		x_j = S \left( \cos\left( \frac{2j \pi}{S} \right) -1 \right), \quad j = 0, \dots, 
		\begin{cases}
			S/2     & S \text{ even } \\
			(S-1)/2 & S \text{ odd }
		\end{cases}.
	\end{equation}
	We give two motivations for considering the Chebyshev extreme points in connection with the pseudo-extrema.
	First, as outlined in \ref{sec:PECheby} the pseudo-extrema of (shifted) Chebyshev polynomials are trivially given by the (shifted) positive Chebyshev extreme points.
	This is of significance since the first order optimal stability polynomial for parabolic spectra is precisely given by the shifted Chebyshev polynomial of first kind.
	Second, we recall that the Chebyshev extreme points $x_j$ can be perceived as the real part of the set of points $z_j$ which partition the 
	circle into segments with equal arc length \cite{trefethen2019approximation}.
	In that case, the imaginary part $y_j$ of $z_j$ can be computed as 
	\begin{equation}
		\label{eq:ImagPartChebyExtremePoints}
		y_j = \pm \sqrt{ S^2 - (x_j + S)^2} = \pm S \sqrt{ 1 - \cos^2\left( \frac{2j \pi}{S} \right)} = \pm S \sin \left( \frac{2j \pi}{S} \right),
	\end{equation}
	since $\sin\left( \frac{2j \pi}{S} \right) \geq 0, j = 0, \dots, S/2$.
	Consequently, the points $z_j = x_j \pm i y_j$ with $x_j, y_j$ defined as above lie on the stability boundary $\partial \mathcal S$ and are thus a valid candidate for the pseudo-extrema. 
	In fact, one can readily show that these points are indeed the pseudo-extrema of $P^\text{Disk}_{S, 1}(z)$, as stated in the following theorem.
	\begin{thm}
		\label{thm:thm1}
		The $S-1$ pseudo-extrema of $P^\text{Disk}_{S, 1}(z)$ are given by the positive Chebyshev extreme points with $x_j \neq 0$ of $T_S(1+z/S)$ with projection onto the circle with radius $S$ centered at $(-S, 0)$.
	\end{thm}
	\begin{pf}
		\label{thm:Proof1}
		We compute 
		\begin{align}
			P^\text{Disk}_{S, 1}(\widetilde{ r}_j) =  \left(1 + \frac{\widetilde{ r}_j}{S}\right)^S &= \left(1 + \frac{x_j \pm i y_j}{S} \right)^S \\
			&= \left(1 + \frac{S \left( \cos\left( \frac{2j \pi}{S} \right) -1 \right) \pm i S \sin \left( \frac{2j \pi}{S} \right)}{S} \right)^S\\
			&= \left(\cos\left( \frac{2j \pi}{S} \right) \pm i \sin \left( \frac{2j \pi}{S} \right) \right)^S \\
			&= \exp( \pm \, i \, 2j \, \pi ) \overset{j \in \mathbb N}{=} 1.
		\end{align}
		As a consequence, we have that 
		\begin{equation}
			\widetilde{P}_{S-1}(\widetilde{ r}_j; \widetilde{\boldsymbol r}) = \frac{P^\text{Disk}_{S,1}(\widetilde{\boldsymbol r}_j; \widetilde{\boldsymbol r}) - 1}{\widetilde{ r}_j} = 0 \quad j = 1, \dots 
			\begin{cases}
				S/2     & S \text{ even } \\
				(S-1)/2 & S \text{ odd }
			\end{cases}
		\end{equation}
		where we use the fact that $P_{S, 1}(0, \widetilde{\boldsymbol r}) = 1$ which excludes $j=0$ from the set of pseudo-extrema.
	\end{pf}
	\begin{rmk}
		Note that for even stability polynomials we have one purely real pseudo extremum $\widetilde{r}_0 = -S$ and $S-2$ complex conjugated pseudo-extrema.
		For odd stability polynomials, all pseudo-extrema are complex conjugated.
		This is in accordance with the motivating observation \cref{eq:PseudoExtremaSubsetStabBnd} where we established that the pseudo-extrema lie on the boundary of the region of absolute stability $\partial \mathcal S$.
		We also note that the complex-conjugated pseudo-extrema of $P^\text{Disk}_{S, 1}(z)$ correspond to the real pseudo-extrema with root-multiplicity two of the Chebyshev polynomials $T_S(1+z/S^2)$, see \ref{sec:PECheby}.
	\end{rmk}
	\subsubsection{Second Order Consistent Stability Polynomial}
	For disks/circular spectra there is also an explicit formula for the optimal second order accurate polynomial known \cite{owren1990some}.
	The polynomial
	\begin{equation}
		\label{eq:StabPolyCircleSecondOrderAccurate}
		P^\text{Disk}_{S, 2}(z) = \frac{S-1}{S} \left( 1 + \frac{z}{S-1}\right)^S + \frac{1}{S}
	\end{equation}
	contains the maximum disk $D_{S-1}$ in its region of absolute stability \cite{owren1990some}.
	As a side note, we realize that $P_{S, 2}(z)$ is identical to the second order accurate Runge-Kutta Chebyshev method \cite{verwer1990convergence, van1980internal} with infinite damping \cite{verwer2004rkc}.
	Furthermore, as mentioned in \cite{verwer2004rkc} $P_{S, 2}(z)$ can be written as a sequence of forward Euler steps and belongs consequently to the family of \ac{tvd}/\ac{ssp} integrators \cite{shu1988efficient, gottlieb2001strong} with optimal \ac{ssp} coefficient \cite{doi:10.1142/7498}.
	\begin{thm}
		The $S-1$ pseudo-extrema of $P^\text{Disk}_{S, 2}$ are given by the positive Chebyshev extreme points $x_j \neq 0$ of $T_S\big(1+z/(S-1)\big)$ when projected onto the circle with radius $S-1$ centered at $\big(-(S-1), 0 \big)$.
	\end{thm}
	\begin{pf}
		Analogous to 
		proof of Theorem \ref{thm:thm1}.
	\end{pf}
	\begin{rmk}
		As for the first order accurate stability polynomial, we have  for even $S$ one real pseudo-extremum $\widetilde{r}_0 = -(S-1)$ and $S-2$ complex conjugated pseudo-extrema and for odd $S$ a set of $(S-1)/2$ exclusively complex conjugated pseudo-extrema.
	\end{rmk}

	The significance of these results is that the pseudo-extrema are distributed with exact equal arc length on the hull of the spectrum, and consequently with approximately equal arc length on the stability boundary for a finite number of stages $S$.
	In fact, in the limit $S\to \infty$ the stability boundaries converge to the circle and thus the pseudo-extrema are asymptotically distributed with exact equal arc length on th stability boundary.
	
\begin{figure}[!t]
		\centering
		\subfloat[Pseudo-Extrema $\widetilde{r}_j$ and boundary of region of absolute stability $\partial \mathcal S$ of $P^\text{Disk}_{12, 2}(z)$ alongside positive Chebyshev extreme points of $T_{12}(1+z/S)$.]{
			\label{fig:CircleEven}
			\centering
			\includegraphics[width=0.46\textwidth]{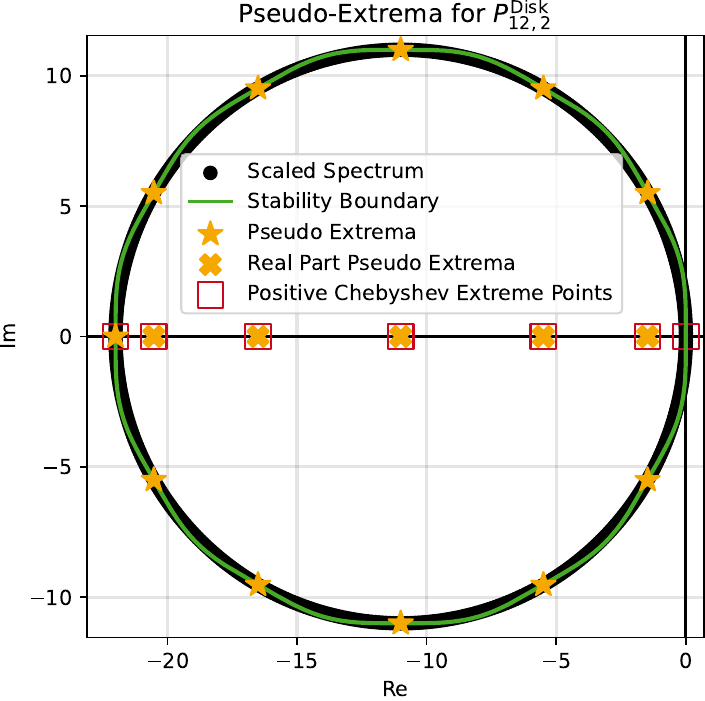}
		}
		\hfill
		\subfloat[Pseudo-Extrema $\widetilde{r}_j$ and boundary of region of absolute stability $\partial \mathcal S$ of $P^\text{Disk}_{13, 2}(z)$ alongside positive Chebyshev extreme points of $T_{13}(1+z/S)$.]{
			\label{fig:CircleOdd}
			\centering
			\includegraphics[width=0.46\textwidth]{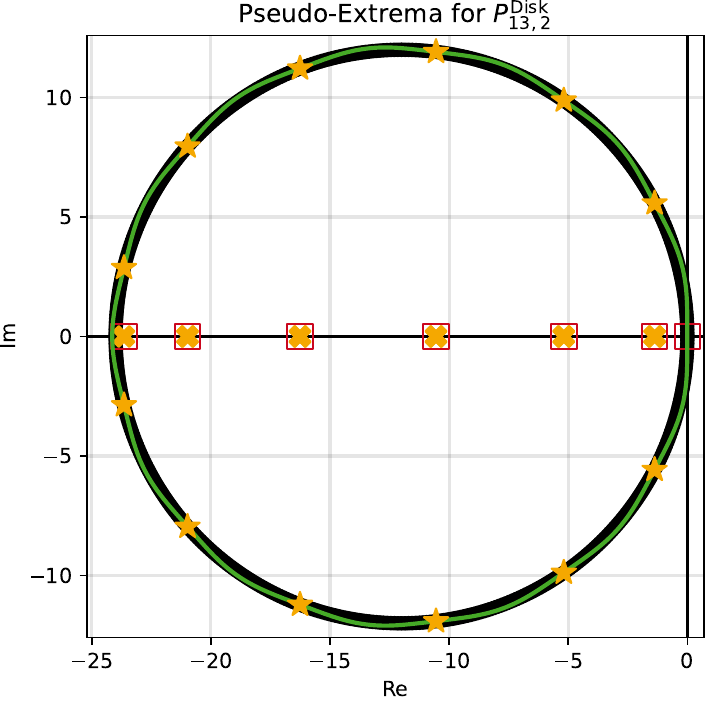}
		}
		\caption[Pseudo-Extrema for Circular Spectrum]{Optimal second order accurate stability polynomials $P_{S, 2}(z)$ given by \cref{eq:StabPolyCircleSecondOrderAccurate} for even \subref{fig:CircleEven} and odd degree \subref{fig:CircleOdd} polynomial.
			The spectrum corresponds to the canonical first order Finite Volume Upwind/Godunov discretization \cite{godunov:hal-01620642} of the advection equation $u_t + u_x = 0$ on the periodic $[-1, 1]$ domain discretized with $500$ cells.
			In the even degree case \subref{fig:CircleEven}, the lower degree polynomial is odd and thus there is one pure real pseudo extremum at the left end of the spectrum.
			For odd polynomial degrees \subref{fig:CircleOdd}, the lower degree polynomial is of even degree and we have only complex-conjugated pseudo-extrema.
			Note also that the segment crossing $0$ is twice the length of the others, which follows from the fact that $P_{S, p}(0) \equiv 1 \: \forall \: S, p \geq 1$, cf. \cref{eq:DefinitionLowerDegree}.
		}
		\label{fig:Circles}
	\end{figure}

	Especially higher degree stability polynomials have enough flexibility to fully adapt the stability boundary to the envelope of the spectrum.
	This is the case for even and odd polynomial degrees $S$ as well as first and second order accurate polynomials.
	To illustrate this, we present \cref{fig:Circles} where the pseudo-extrema of second order accurate polynomials $P^\text{Disk}_{S, 2}$ with degrees $S=12$ and $S=13$ are shown.
	\subsection{Central Assumption}
	\label{sec:CentralAssumption}
	Based on the findings for the circle we make a central assumption regarding the pseudo-extrema of strictly convex spectra, i.e., where the eigenvalues define a strictly convex set of points in the complex plane.
	To establish an intuition for this type of spectra we present four exemplarily representatives falling into this category in \cref{fig:ConvexSpectra}.
	We assume that the pseudo-extrema of 
	optimal stability polynomials are distributed with approximately equal arc length on the convex hull 
	of strictly convex spectra.
	\subsubsection{Strictly Convex Spectra}
\begin{figure}[!t]
		\centering
		\subfloat[{Burgers' Equation corresponding to the semidiscretization defined by item \ref{item:Burgers}.}]{
			\label{fig:ConvexSpectraBurgers}
			\centering
			\includegraphics[width=0.45\textwidth]{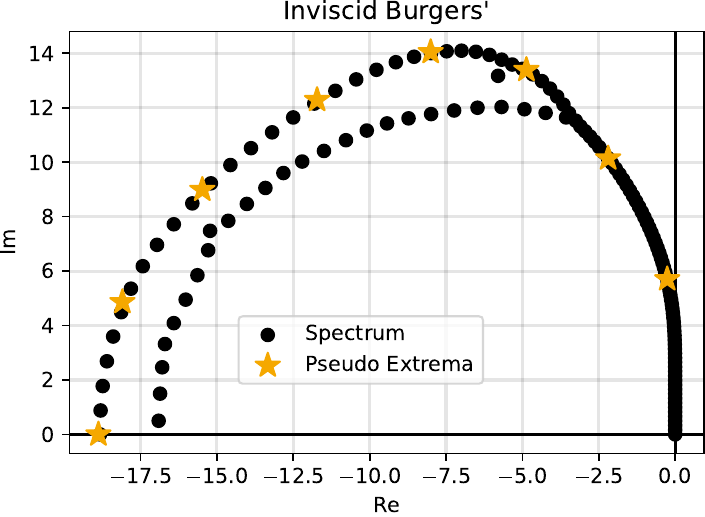}
		}
		\hfill
		\subfloat[{1D Shallow Water Equations corresponding to the semidiscretization defined by item \ref{item:ShallowWater}.}]{
			\label{fig:ConvexSpectraShallowWater}
			\centering
			\includegraphics[width=0.45\textwidth]{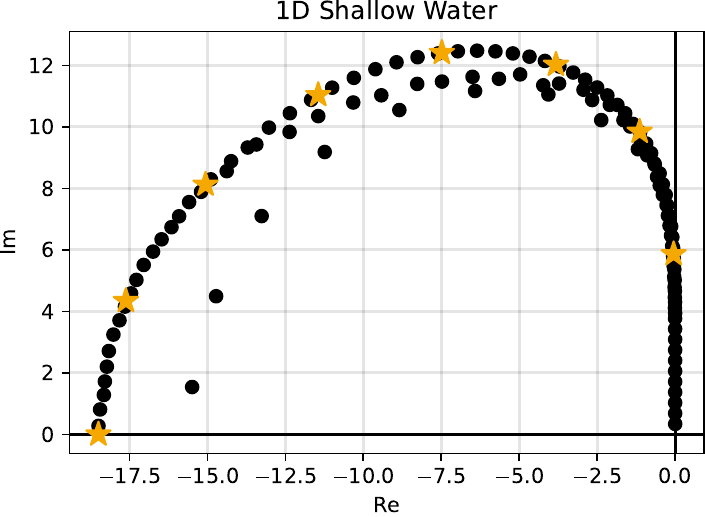}
		}
		\\
		\subfloat[{1D Compressible Ideal MHD corresponding to the semidiscretization defined by item \ref{item:MHD}.}]{
			\label{fig:ConvexSpectraMHD}
			\centering
			\includegraphics[width=0.45\textwidth]{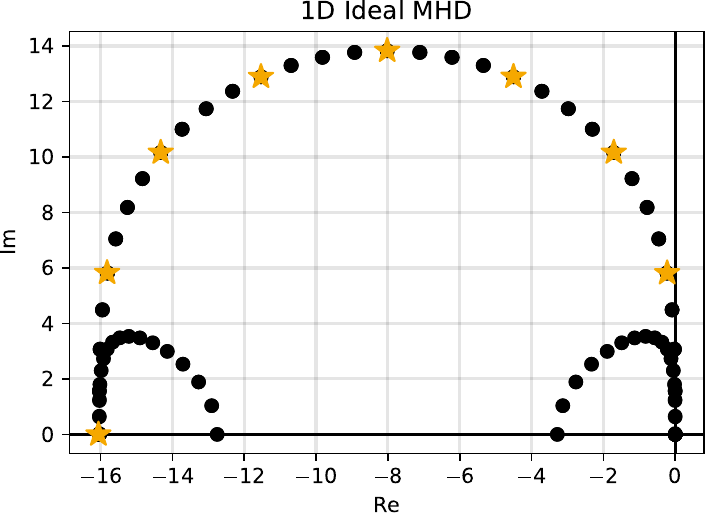}
		}
		\hfill
		\subfloat[{2D Compressible Euler equations corresponding to the semidiscretization defined by item \ref{item:2D_CEE}.}]{
			\label{fig:ConvexSpectra2D_CEE}
			\centering
			\includegraphics[width=0.45\textwidth]{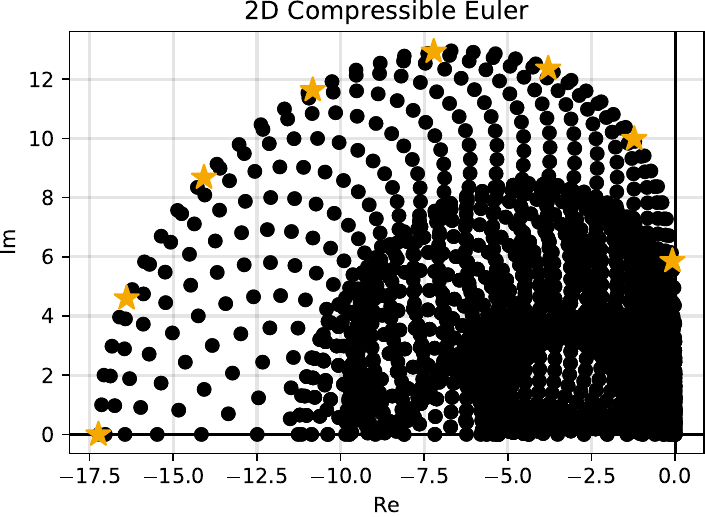}
		}
		\caption[Collection of Convex Spectra]{Collection of strictly convex spectra of nonlinear hyperbolic \ac{pdes} and optimized pseudo-extrema $\widetilde{r}_j$.
			The spectra are scaled with the optimal timestep $\Delta t_{16, 2}$ in each case.
		}
		\label{fig:ConvexSpectra}
	\end{figure}

	To motivate this assumption extending the findings for proven optimal spectra of the circle, we display the pseudo-extrema of optimized $P_{16,2}(z)$ polynomials for different spectra in \cref{fig:ConvexSpectra}.
	The optimal stability polynomials and optimal timesteps $\Delta t_{16,2}$ have been computed using the approach developed in \cite{ketcheson2013optimal} and the pseudo-extrema are computed numerically.
	
	We present the spectra of four nonlinear hyperbolic \ac{pdes}
	\begin{enumerate}
		\item Burgers' equation 
		\item 1D Shallow Water equations
		\item 1D Ideal Compressible \ac{Mhd} equations 
		\item 2D Compressible Euler equations 
	\end{enumerate}
	when discretized with the \ac{dgsem} \cite{kopriva2010quadrature, kopriva2009implementing} using \texttt{Trixi.jl} \cite{trixi1, trixi2, trixi3}.
	We use local polynomials of degrees $1, 2, 3$ and a range of numerical fluxes and initial conditions, corresponding to both smooth and discontinuous solutions.
	For the sake of reproducibility we give the precise setups of the semidiscretizations corresponding to the displayed spectra in \ref{app:SpectraSetup}.
	
	As seen from \cref{fig:ConvexSpectra} the pseudo-extrema are distributed with approximately equal arc length on the spectrum enclosing curve.
	Based on this observation we will initialize the pseudo-extrema such that they are distributed with equal arc length on the spectrum enclosing curve, which is assumed to be very close to the optimal stability boundary of higher degree polynomials.
	The deviations from the equal arc-length initialization are then courtesy of the optimizer.

	The significance of the strictly convex spectra stems from the observation that the optimal stability boundary $\partial \mathcal S^\star$ can be, for a sufficiently high polynomial degree $S$, reasonably well approximated with the convex hull of the eigenvalues.
	This implies, in turn, that we can use the convex hull of the spectrum as an Ansatz for the initial placement of the pseudo-extrema $\widetilde{\boldsymbol r}$.
	Generalizing the proven results for the circle, the pseudo-extrema are then placed with equal arc length on the convex hull of the spectrum.
	\section{Detailed Formulation of the Optimization Problem}
	\label{sec:Formulation}
	In this section we discuss the required details to successfully optimize a stability polynomial when parametrized in pseudo-extrema \cref{eq:OptProbPE}.
	\subsection{Constructing Complex-Conjugated Pseudo-Extrema}
	Since the lower-degree polynomial $\widetilde{P}_{S-1}$ has real coefficients, one has to ensure that the pseudo-extrema $\widetilde{\boldsymbol r}$ are complex-conjugated.
	In principle, this is a difficult task since throughout the optimization it might happen that a previously real pseudo-extremum is moved away from the real axis into the complex plane.
	Consequently, another real pseudo-extremum has to be set as the corresponding complex-conjugate, effectively changing the number of free optimization variables by introducing an additional constraint.
	To avoid this switch-like behaviour we set a-priori a number of real and complex-conjugated pseudo-extrema.
	Based on the previous results, we restrict ourselves for even degree stability polynomials one real and $S-2$ complex conjugated pseudo-extrema and for odd stability polynomials only $S-1$ complex conjugated ones.
	For strictly convex spectra this is always found to be true, cf. \cref{fig:ConvexSpectra}.
	\subsection{Restriction to Second Quadrant}
	Due to symmetry around the real axis it suffices to consider either the second or third quadrant in the complex plane only.
	In this work, we choose the second quadrant, which implies that we consider only complex numbers with non-negative imaginary part.
	This implies that for real Jacobians \cref{eq:JacobianSemiDisc} with complex-conjugated eigenvalues $\lambda_j = \overline{\lambda_{j+1}}$ the stability requirement $\vert P_{S}(\Delta t \lambda^{(m)}) \vert \leq 1$ needs only to be checked for the $m = 1, \dots, \widetilde{M}$ eigenvalues with positive imaginary part, thereby also reducing the number of constraints.
	\subsection{Scaling of the Optimal Timestep}
	We briefly comment on the scaling of the optimal timestep $\Delta t^\star$ with the degree of the stability polynomial $S$.
	It is well-known that the maximum admissible timestep scales for parabolic spectra covering the real-axis quadratic in $S$ \cite{abdulle2015explicit, van1996development}.
	In contrast, the one-sided imaginary axis inclusion of polynomials with real coefficients is bound by $S-1$ as shown in \cite{vichnevetsky1983new}.
	In \cite{kinnmark1984one2} the tighter maximum imaginary stability limit $\sqrt{(S-1)^2 -1}$ was conjectured and proven for special cases in \cite{kinnmark1987principle}.
	The proven optimal stability polynomials for the disk scale both linearly in $S$, see \cite{jeltsch1978largest, owren1990some}.
	As a consequence, we consider for general spectra an asymptotically linear scaling
	\begin{equation}
		\label{eq:TimestepLinearScaling}
		\Delta t_{S,p} = \frac{S}{S_\text{Ref}} \Delta t_{S_\text{Ref},p},
	\end{equation}
	which is observed for instance in \cite{ketcheson2013optimal, VERMEIRE2019465, Nasab2022Optimal, al2021optimized, hedayati2021optimal}.
	It should be mentioned that for $\mathcal{O}(1)$ better than linear scalings are possible, as for these degrees qualitative changes in the shape of the stability boundary still occur. 
	Once the stability boundary is closely adapted to the spectrum, the linear scaling is recovered.
	
	It is immediately clear that once the timestep increases only linearly with more stages $S$ there is no additional efficiency gained.
	In particular, we observe for the methods with a large number of stages additional complications due to internal stability, see \cref{subsec:InternalStabiliy}.
	Nevertheless, many stage stability polynomials are of great value in the context of multirate partitioned Runge-Kutta.
	As mentioned before, the \ac{PERK} schemes \cite{VERMEIRE2019465, HEDAYATINASAB2022111470, VERMEIRE2023112159} would benefit from high degree stability polynomials which enable the efficient integration of systems even in the presence of locally restricted CFL numbers.
	\subsection{Constraints for Higher Order}
	For higher order $p$-consistent methods $p-1$ equality constraints have to be added to the stability inequality constraints.
	To obtain second order consistent methods, for instance, the pseudo-extrema have to satisfy
	\begin{equation}
		\label{eq:SecOrdConsConstr}
		\frac{1}{2} \overset{!}{=} - \sum_{j=1}^{S-1} \frac{1}{\widetilde{r}_j}
	\end{equation}
	which follows from Vieta's formulas which establish a link between the roots and coefficients of a polynomial.
	More generally, for a stability polynomial that matches the first $k=0, \dots, p$ coefficients of the exponential $1/(j!)$ we have the additional constraints
	\begin{equation}
		\label{eq:GeneralFormOrderConstraint}
		\frac{1}{k!} \overset{!}{=} (-1)^{k-1}\sum_{\substack{1=j_1, j_2 \dots j_{k-1}\\j_1 \neq j_2 \neq \cdots \neq j_{k-1}}}^{S-1} \frac{1}{\widetilde{\boldsymbol r}_{j_1} \widetilde{\boldsymbol r}_{j_2} \dots \widetilde{\boldsymbol r}_{j_{k-1}}}, \quad k = 2, \dots, p.
	\end{equation}
	In this work we focus on second order accurate polynomials since for these the linear order constraints \eqref{eq:GeneralFormOrderConstraint} imply second order convergence also for nonlinear system, without any additional constraints placed on the coefficients of the method \cite{HairerWanner1, wanner1996solving}.
	Nevertheless, we also constructed third order accurate stability polynomials for some selected cases.
	\subsection{Focus on Even Degree Stability Polynomials}
	\label{subsec:EvenDegree}
	From now on we focus on even degree polynomials since they lead to a slightly simpler formulation of the optimization problem.
	As we are interested in the many stage case, focusing on either even or odd degree case is no severe restriction.
	
\begin{figure}[!t]
		\centering
		\includegraphics[width=0.45\textwidth]{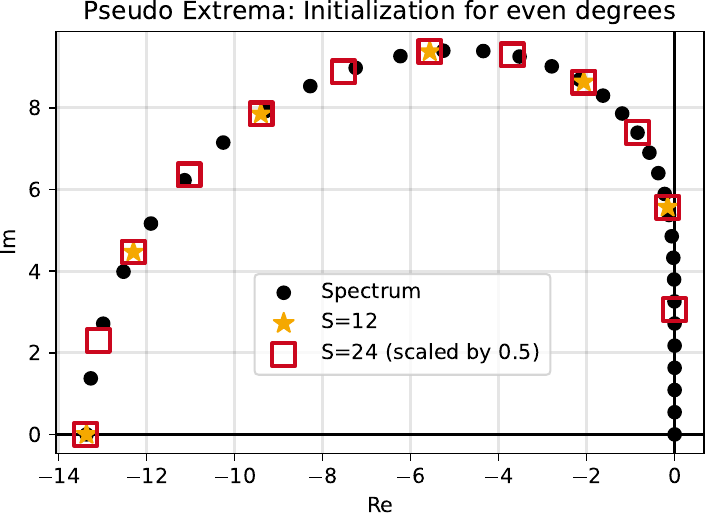}
		\caption[Re-Usability of Pseudo-Extrema of Optimal Even Degree Stability Polynomials.]{Scaled spectrum $\sigma_{12, 2}$ of the 1D linear advection equation $u_t +u_x =0$ discretized through the \ac{dgsem} on $[0, 1]$ using $16$ cells/elements with DG polynomial degree $3$ and Rusanov/Local Lax-Friedrichs flux.
		The pseudo-extrema of the $12$ and $24$ degree second order accurate polynomial are computed.
		To highlight that the relative positions of every second pseudo extremum $\tilde{r}_{2j-1}^{2S}$ agree with $\tilde{r}_j^S$, the former are scaled by $0.5$ due to the linear scaling of the timestep according to \cref{eq:TimestepLinearScaling}.
		}
		\label{fig:DoublingEven}
	\end{figure}

	Even degree stability polynomials come with two main advantages.
	First, we observe for all considered spectra (not only strictly convex) that the left endpoint of the spectrum, in our cases a real eigenvalue, is always selected as one pseudo-extremum, allowing for an even improved initialization.
	For strictly convex spectra, this is the only real pseudo-extremum, cf. \cref{fig:CircleEven}, \cref{fig:ConvexSpectra}, and \cref{fig:DoublingEven}.
	The remaining pseudo-extrema can then be initialized with equal arc length starting from the left endpoint of the spectrum.
	
	Second, under the assumption motivated in \cref{sec:CentralAssumption} and the linear scaling of the timestep one can readily re-use the optimal pseudo-extrema of the $S/2$ degree polynomial as an advanced initialization for the $S$ degree polynomial.
	In particular, every second pseudo-extremum of the $S$ stability polynomial is essentially known, cf. \cref{fig:DoublingEven}.
	\subsection{Construction of the Convex Hull Interpolant}
	For strictly convex spectra the convex hull defines in the second quadrant actually a strictly concave (strictly convex in the third quadrant) function, which allows to uniquely map the real part of the pseduo-extrema to some corresponding imaginary part.
	To relax the very steep section near the imaginary axis we add the origin $(0,0)$ to the hull which helps to reduce the slope in that particular segment.
	This allows to carry out the optimization actually in terms of the real parts of the pseudo-extrema only, cutting the number of optimization variables in half.
	
	Given the (strict) convex hull of the spectrum defined through a set of eigenvalues $\mu^{(j)} \in \sigma^\text{Hull}$ and the expected maximum timestep $\Delta t_{S,p}$, we define the scaled spectrum 
	\begin{align}
		\sigma_{S,p}^\text{Hull} &\coloneqq \Big\{ \Delta t_{S,p} \cdot \lambda, \: \lambda \in \sigma^\text{Hull} \Big\} \\ 
		\mu_{S,p}^{(j)} &\coloneqq \Delta t_{S,p} \cdot \mu^{(j)}\:.
	\end{align}
	Based on this we construct a piece-wise linear interpolant $I$ as
	\begin{align}
		\begin{split}
			\label{eq:PieceWiseLinearIntPol}
			I\left(x; \sigma^\text{Hull}_{S,p} \right) \coloneqq \: & \text{Im}\left(\mu_{S,p}^{(j)} \right) + \frac{\text{Im}\left(\mu_{S,p}^{(j+1)}\right) -\text{Im}\left(\mu_{S,p}^{(j)}\right)}{\text{Re}\left(\mu^{(j+1)}_{S,p}\right) -\text{Re}\left(\mu_{S,p}^{(j)}\right)} \left(x - \text{Re}\left(\mu_{S,p}^{(j)}\right)\right), \\
			&\text{Re}\left(\mu_{S,p}^{(j)}\right) \leq x < \text{Re}\left(\mu_{S,p}^{(j+1)}\right).
		\end{split}
	\end{align}
	While quadratic or the classic cubic splines would result in a continuously differentiable interpolant $I$, they lead to numerical artifacts in the very steep sections near the imaginary axis that are a common feature of the \ac{dg} spectra.
	For these parts it is observed that linear splines capture the boundary of the stability domain $\partial \mathcal{S}$ best and should be used there.
	For the less steep regions quadratic and cubic splines were tested, but the increase in approximation quality compared to linear splines was for the majority of cases too small to justify the increased computational costs.
	
	Global interpolation techniques were also tested but familiar problems like Runge's phenomenon and large over/undershoots in the intervals between the eigenvalues are observed for a variety of approaches, including Bernstein polynomials or mapping to Chebyshev points \cite{DEMARCHI2020112347}.
	
	It should be mentioned that the spectrum close to the imaginary axis becomes arbitrarily steep and thus the interpolation becomes less reliable in this part.
	This becomes an issue for stability polynomials of degrees $\mathcal O(100)$ since then pseudo-extrema are placed in that region.
	To solve this problem, one can work with a spectrum enclosing curve
	$\gamma(\tau): \mathbb R \to \mathbb C$ instead of a spectrum enclosing function $I:\mathbb R \to \mathbb R$.
	By identifying the eigenvalues on the hull $\mu^{(j)} \in \sigma^{\text{Hull}}$ with an arc length parameter $0 \leq \tau^{(j)} \leq 1$ the imaginary part can be computed as
	\begin{align}
		\begin{split}
			\label{eq:PieceWiseLinearIntPolAlphaShapeImag}
			I_\text{Im}\left(\tau; \sigma^{\text{Hull}}_{S,p} \right) = \:& \text{Im}\left(\mu_{S,p}^{(j)} \right) + \frac{\text{Im}\left(\mu_{S,p}^{(j+1)}\right) -\text{Im}\left(\mu_{S,p}^{(j)}\right)}{\tau^{(j+1)}_{S,p} -\tau_{S,p}^{(j)}} \left(\tau - \tau_{S,p}^{(j)} \right), \\
			&\tau_{S,p}^{(j)} \leq \tau < \tau_{S,p}^{(j+1)}
		\end{split}
	\end{align}
	which removes the vanishingly small difference in real parts from the denominator.
	This, however, comes at the expense that we have to do the interpolation \cref{eq:PieceWiseLinearIntPolAlphaShapeImag} also for the real part of the pseudo-extrema, thus doubling the computational work load.
	As discussed in \cref{sec:ResultsStrictlyConvexSpectra}, methods with more than hundred stages introduce additional complications related to internal stability.
	Consequently, we tailor the optimization approach around the more robust $S \sim \mathcal O(10)$ case and use the real parts of the pseudo-extrema as the optimization variables.
	\section{Implementation}
	\label{sec:Implementation}
	In this section, we bundle the aforementioned details into the precise formulation of the solved problems.
	\subsection{Optimization Problems}
	Equipped with the interpolant $I(x)$ the optimization of lower-degree polynomial can now be efficiently conducted.
	The optimization problem \cref{eq:OptProb} is first conducted in the real parts $x_j = \text{Re}(\widetilde{r}_j)$ only, i.e.,
	\begin{subequations}
		\label{eq:OptProbx}
		\begin{align}
			\label{eq:OptProbNonConvexSpectraObjective}
			&\max \Delta t \text{ over }
			\boldsymbol x \in \left[ \underset{m}{\min}\left\{\text{Re}\left(\lambda_{S,p}^{(m)}\right) \right\}, 0\right]^{N} \text{ such that }  \\
			\label{eq:OptProbNonConvexSpectraConstr}
			&\left \vert 1 + \left(\Delta t \lambda^{(m)}\right)  \widetilde{P}_{S-1}\left(\Delta t \lambda^{(m)}; \boldsymbol x + iI\Big(\boldsymbol x; \sigma^\text{Hull}_{S,p}\Big) \right) \right \vert \leq 1 \quad  m = 1 , \dots, \widetilde{M}.
		\end{align}
	\end{subequations}
	Note that the number of optimization variables $N$ equals for even stability polynomials $1 + \frac{S-2}{2} = \frac{S}{2}$ since we expect one real pseudo extremum and $S-2$ complex conjugated ones.
	
	Due to the sensitive dependence of the stability and order constraints on the pseudo-extrema we conduct a second optimization run where we allow for small corrections $y$ in the imaginary parts of the pseudo-extrema:
	\begin{equation}
		\boldsymbol y \in \left[-\eps \: \underset{m}{\max}\left\{\text{Im}\left(\lambda_{S,p}^{(m)}\right) \right\}, \eps \: \underset{m}{\max}\left\{\text{Im}\left(\lambda_{S,p}^{(m)} \right) \right\} \right]^{N}.
	\end{equation}
	Here, $\eps$ is the sole hyperparameter needed in our approach which could for all problems be set to $\mathcal O \left(10^{-2}\right)$ with default choice $0.02$.
	The second stage optimization problem reads then
	\begin{subequations}
		\label{eq:OptProbxy}
		\begin{align}
			&\max \Delta t \text{ over } \begin{pmatrix} \boldsymbol x \\ \boldsymbol y \end{pmatrix} \in \begin{pmatrix} \left[ \underset{m}{\min}\left\{\text{Re}\left(\lambda_{S,p}^{(m)}\right) \right\}, 0\right]^{N} \\
				\left[-\eps \: \underset{m}{\max}\left\{\text{Im}\left(\lambda_{S,p}^{(m)}\right) \right\}, \eps \: \underset{m}{\max}\left\{\text{Im}\left(\lambda_{S,p}^{(m)} \right) \right\} \right]^{N} \end{pmatrix} \\ &\text{such that } \nonumber \\
			\label{eq:StabConstraintFeasxy}
			&\left \vert 1 + \left(\Delta t \lambda^{(m)}\right)  \widetilde{P}_{S-1}\left(\Delta t \lambda^{(m)}; \boldsymbol x + i \left[I\Big(\boldsymbol x; \sigma^\text{Hull}_{S,p}\Big) + \boldsymbol y \right] \right) \right \vert \leq 1 \quad  m = 1 , \dots, \widetilde{M}.
		\end{align}
	\end{subequations}
	The entire algorithm is given in \cref{alg:RootOpt}.
	In principle, it is also possible to conduct the optimization in real parts and imaginary corrections \cref{eq:OptProbxy} only, i.e., without solving \cref{eq:OptProbx} first.
	In all considered cases, however, optimizing the real parts of the pseudo-extrema first led not only to a more robust, but also overall faster optimization process.
	\begin{algorithm2e}[!ht]
		\SetKw{KwAnd}{and}
		\nl Perform Algorithm proposed in \cite{ketcheson2013optimal} to obtain $\Delta t_{S_\text{Ref},p}$ for certain $S_\text{Ref}, p$\;
		\nl$\Delta t_{S,p} \leftarrow \Delta t_{S_\text{Ref},p} \frac{S}{S_\text{Ref}}$\;
		\nl Scale spectrum with expected timestep: $\sigma^{(S,p)} \coloneqq \{ \lambda^{(m)} \cdot \Delta t_{S,p} \}_{m= 1, \dots, \widetilde{M}}$\;
		\nl Compute convex hull $\sigma^\text{Hull}_{S,p} \leftarrow \gamma_0\left(\sigma^{(S,p)}\right)$\;
		\nl Compute Interpolant $I\left(x; \sigma^\text{Hull}_{S,p}\right)$\;
		\nl \If{$S, S/2 $ are even \KwAnd results $\boldsymbol x^{(S/2)}$ for the $S/2$ degree polynomial exist}{
			\nl Set every second entry in $\boldsymbol x_0^{(S)}$ to $2 \boldsymbol x^{(S/2)}$ \;
			\nl Set remaining entries such that $\left( x_{0,j}, I\left( x_{0,j}; \sigma^\text{Hull}_{S,p}\right)\right), j = 2, 4, \dots, N$ are equally distributed on $I\left(x; \sigma^\text{Hull}_{S,p}\right)$\;
		}
		\nl\Else{
			\nl Compute initial guess $\boldsymbol x_0$ such that $\left( x_{0,j}, I\left(x; \sigma^\text{Hull}_{S,p}\right) \right), j = 1, \dots, N$ are equally distributed on $I\left(x; \sigma^\text{Hull}_{S,p}\right)$\;
		}
		\nl Initialize timestep with maximum possible value $\Delta t_0 \leftarrow \Delta t_{S,p}^\star$\;\label{line:11}
		\nl Find best candidate $\Delta t^\star, \boldsymbol x^*$ of optimization problem \eqref{eq:OptProbx}\;\label{line:12}
		\nl Update initial guess $\Delta t_0 \leftarrow \Delta t^\star, \boldsymbol x_0 \leftarrow \boldsymbol x^*, \boldsymbol y_0 = \boldsymbol 0$ \;
		\nl $(\boldsymbol x^\star, \boldsymbol y^\star ) \leftarrow $ Solve second stage optimization problem \eqref{eq:OptProbxy}\;\label{line:14}
		\nl\Return{$\widetilde{\boldsymbol r} \leftarrow \boldsymbol x^\star \pm i\left[I\left(\boldsymbol x^\star ; \sigma^\mathrm{Hull}_{S,p} \right) + \boldsymbol y^\star\right]$}\;
		\nl\Return{$\Delta t_{S,p}^\star \leftarrow \Delta t^\star$}\;\label{line:16}
		\caption{Optimization problem in pseudo-extrema.}
		\label{alg:RootOpt}
	\end{algorithm2e}
	\subsection{Feasibility Problems}
	Under assumption that we can indeed realize the optimal timestep $\Delta t_{S,p}$ one can formulate the optimization problems \cref{eq:OptProbx}, \cref{eq:OptProbxy} actually as feasibility problems by fixing $\Delta t = \Delta t_{S,p}$ and search for $\boldsymbol x, \boldsymbol y$ satisfying the stability and order constraints.
	The feasibility problem corresponding to \cref{eq:OptProbx} reads then
	\begin{subequations}
		\label{eq:FeasProbx}
		\begin{align}
			\label{eq:FeasProbNonConvexSpectraObjective}
			&\text{Find }
			\boldsymbol x \in \left[ \underset{m}{\min}\left\{\text{Re}\left(\lambda_{S,p}^{(m)}\right) \right\}, 0\right]^{N} \text{ such that }  \\
			\label{eq:FeasProbNonConvexSpectraConstr}
			&\left \vert 1 + \left(\Delta t_{S,p} \lambda^{(m)}\right)  \widetilde{P}_{S-1}\left(\Delta t_{S,p} \lambda^{(m)}; \boldsymbol x + iI\Big(\boldsymbol x; \sigma^\text{Hull}_{S,p}\Big) \right) \right \vert \leq 1 \quad  m = 1 , \dots, \widetilde{M}.
		\end{align}
	\end{subequations}
	and the analogy to \cref{eq:OptProbxy} is given by
	\begin{subequations}
		\label{eq:FeasProbxy}
		\begin{align}
			&\text{Find } \begin{pmatrix} \boldsymbol x \\ \boldsymbol y \end{pmatrix} \in \begin{pmatrix} \left[ \underset{m}{\min}\left\{\text{Re}\left(\lambda_{S,p}^{(m)}\right) \right\}, 0\right]^{N} \\
				\left[-\eps \: \underset{m}{\max}\left\{\text{Im}\left(\lambda_{S,p}^{(m)}\right) \right\}, \eps \: \underset{m}{\max}\left\{\text{Im}\left(\lambda_{S,p}^{(m)} \right) \right\} \right]^{N} \end{pmatrix} \\ &\text{such that } \nonumber \\
			&\left \vert 1 + \left(\Delta t_{S,p} \lambda^{(m)}\right)  \widetilde{P}_{S-1}\left(\Delta t_{S,p} \lambda^{(m)}; \boldsymbol x + i \left[I\Big(\boldsymbol x; \sigma^\text{Hull}_{S,p}\Big) + \boldsymbol y \right] \right) \right \vert \leq 1 \quad  m = 1 , \dots, \widetilde{M}.
		\end{align}
	\end{subequations}
	In practice, the feasibility problems are solved significantly faster than their optimization counterparts.
	Note that no conceptual changes to \cref{alg:RootOpt} are necessary, as the timestep handling in Lines \ref{line:11} and \ref{line:16} can be omitted and \cref{eq:FeasProbx} and \cref{eq:FeasProbxy} are solved in Lines \ref{line:12} and \ref{line:14}, respectively.
	\subsection{Software}
	Both optimization and feasibility problem are solved via 
	\begingroup \spaceskip=\fontdimen2\font plus \fontdimen3\font minus \fontdimen4=\fontdimen7 \ttfamily Ipopt: Interioir Point OPTimizer \endgroup
	\cite{wachter2006implementation}, an optimization software designed for general non-linear problems.
	The herein used required linear solver is \texttt{MUMPS: MUltifrontal Massively Parallel sparse direct Solver} \cite{MUMPS1, MUMPS2} with optional package \texttt{METIS} \cite{karypis1998fast}.
	Being a derivative based optimizer, \texttt{Ipopt} requires the Jacobian matrix and Hessian tensor of the constraints (and objective) to construct the Lagrangian of the problem.
	Here, the derivatives are computed algorithmically via \texttt{dco/c++} \cite{naumann2014dco} which requires only the adaption of some boilerplate code.
	Besides being much more computationally efficient than approximating the derivatives using finite differences, algorithmic differentiation provides exact derivatives (up to machine precision).
	This is especially relevant in this case, since the stability constraints are very sensitive with respect to the pseudo-extrema, i.e., small deviations in $\boldsymbol x$ can make the difference between a stable and an unstable method.
	This was observed when \texttt{Matlab}'s function for general nonlinear optimization, \texttt{fmincon} failed in finding the optimal solution.
	Therein, finite difference approximations of Jacobians and Hessians are used and consequently, the optimization is much less reliable than \texttt{Ipopt} combined with \texttt{dco/c++}.
	
	In terms of complexity the proposed algorithms scales linearly in the number of constraints (stability constraints at eigenvalues and $p-1$ order constraints) $\widetilde{M}$ and quadratic in the number of unknowns (pseudo-extrema with distinct real part) $N$.
	
	It should be mentioned that \texttt{Ipopt} provides many options which can significantly speed up the optimization/feasible point search.
	Most notably, one should specify \texttt{grad\_f\_constant yes} which ensures that the gradient of the objective is only evaluated once and the Hessian not at all.
	An option for a constant objective function as the case for the feasibility problems is not supported, but can actually be implemented by changing only two lines of the \texttt{Ipopt} source code \cite{IPOPTConstantObjective}.
	Furthermore, for many problems the Hessian of the Lagrangian may be successfully approximated using the Limited-memory BFGS method \cite{liu1989limited} which is activated by specifying the option \texttt{hessian\_approximation limited-memory}, resulting in significant speed up.
	In terms of tolerances, we demand satisfaction of the constraints to machine precision by setting \texttt{constr\_viol\_tol 1e-16}.
	
	Furthermore, in most cases it suffices to supply only the convex hull and a small subset of the eigenvalues instead of the entire spectrum.
	This can be seen from e.g. \cref{fig:PE_ConvexSpectraBurgers} and \cref{fig:PE_ConvexSpectra2D_CEE} 
	where the stability constraint at the "interior" eigenvalues is usually satisfied without being explicitly enforced.
	One has to ensure, however, that there is always a larger number of constraints than pseudo-extrema to prevent degenerate cases where the roots of the stability polynomial are placed exactly on the eigenvalues, thus allowing in principle infinitely large timesteps.
	
	The source code is made publicly available on \texttt{GitHub} \cite{OSPREI}.
	\section{Optimal Stability Polynomials for Strictly Convex Spectra}
	\label{sec:OptStabPolysStrictlyConvexSpectra}
	Here, we present many-stage optimal stability polynomials corresponding to first and second order for the spectra shown in \cref{fig:ConvexSpectra}.
	For each spectrum, we compute $\Delta t_{16,p}$ and $P_{16, p}(z)$ using the approach developed in \cite{ketcheson2013optimal}.
	Equipped with the optimal timestep, we supply the expected timestep according to \cref{eq:TimestepLinearScaling} to the feasibility problems \cref{eq:FeasProbx}, \cref{eq:FeasProbxy}.
	For every spectrum, we start with a moderate number of stages ($26-32$) and then double the number of stages to utilize the equal arc length property for even spectra discussed in \cref{subsec:EvenDegree}.
	
\begin{figure}[!t]
		\centering
		\subfloat[{Spectrum of Burgers' Equation computed for the semidiscretization defined by item \ref{item:Burgers}.}]{
			\label{fig:PE_ConvexSpectraBurgers}
			\centering
			\includegraphics[width=0.45\textwidth]{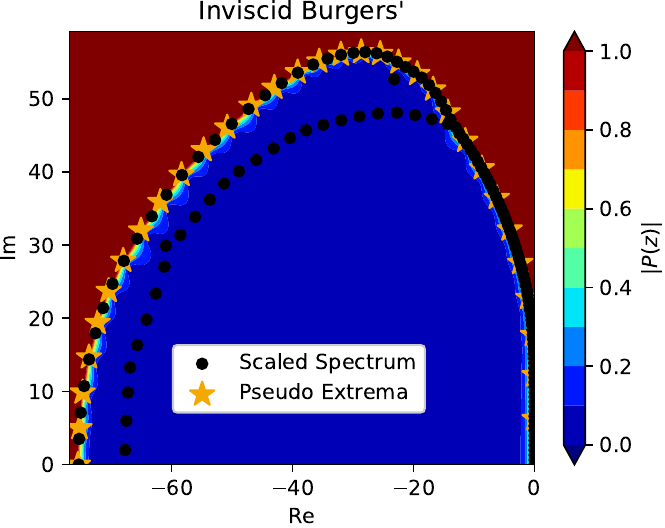}
		}
		\hfill
		\subfloat[{Spectrum of the 1D Shallow Water Equations computed for the semidiscretization defined by item \ref{item:ShallowWater}.}]{
			\label{fig:PE_ConvexSpectraShallowWater}
			\centering
			\includegraphics[width=0.45\textwidth]{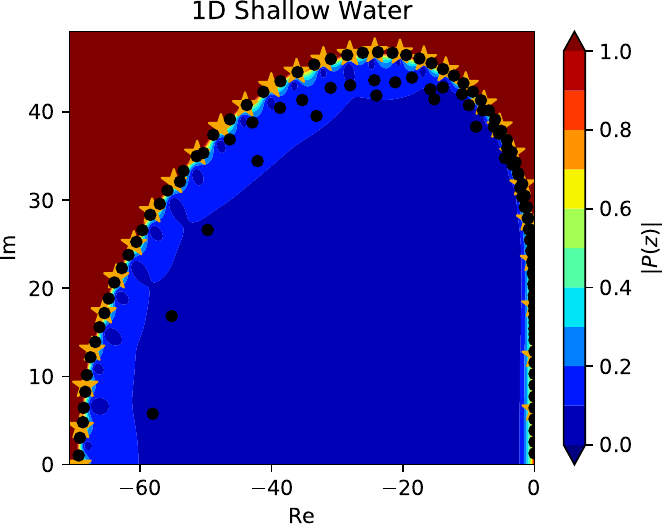}
		}
		\\
		\subfloat[{Spectrum of the 1D Compressible Ideal MHD computed for the semidiscretization defined by item \ref{item:MHD}.}]{
			\label{fig:PE_ConvexSpectraMHD}
			\centering
			\includegraphics[width=0.45\textwidth]{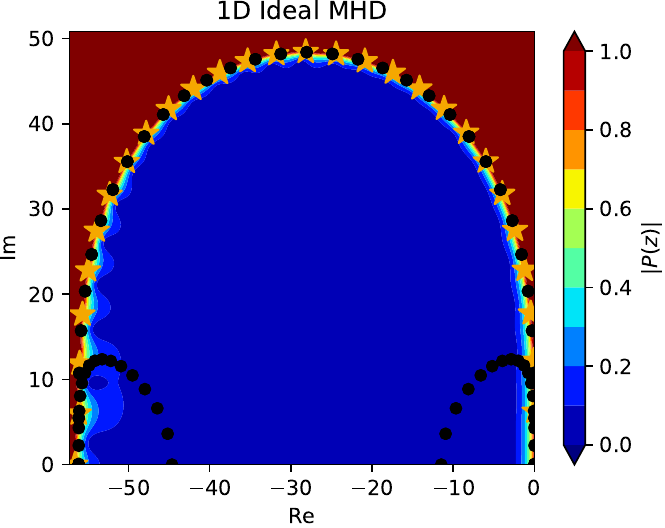}
		}
		\hfill
		\subfloat[{Spectrum of the 2D Compressible Euler equations computed for the semidiscretization defined by item \ref{item:2D_CEE}.}]{
			\label{fig:PE_ConvexSpectra2D_CEE}
			\centering
			\includegraphics[width=0.45\textwidth]{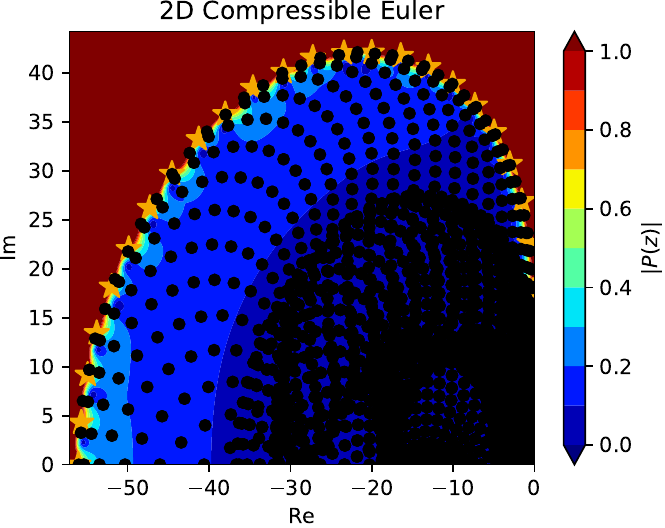}
		}
		\caption[Collection of Convex Spectra]{Collection of strictly convex spectra of nonlinear hyperbolic \ac{pdes} and optimized pseudo-extrema.
			The spectra are scaled with the optimal timestep $\Delta t_{16, 2}$ in each case.
		}
		\label{fig:PE_ConvexSpectra}
	\end{figure}

	The computed pseudo-extrema alongside the constraining list of eigenvalues are provided as supplementary material to this manuscript.
	Optimal stability polynomials are constructed for:
	\begin{enumerate}
		\item Burgers' Equation: We find the optimal stability polynomials of first, second, and third order with degrees $32, 64, 128$ constrained by $\widetilde{M} = 129$ eigenvalues.
		In particular, we make use of the property discussed in \cref{subsec:EvenDegree} by conducting the optimization sequentially, i.e., first performing the optimization for $ S=32$, re-using this results for the initialization of the $S=64$ polynomial and in turn using this results for the final $S=128$ polynomial.
		The pseudo-extrema and contours of the optimal $P_{64,2}(z)$ stability polynomial are displayed in \cref{fig:PE_ConvexSpectraBurgers}.
		
		To highlight the efficiency of our approach, we provide the \texttt{Ipopt} runtimes for the feasibility problem runs in \cref{tab:TimesBurgers}. 
		We emphasize that we do not intend to compare the runtimes of \cref{alg:RootOpt} to the optimization problem from \cite{ketcheson2013optimal} - instead, we would like to stress that equipped with an optimal timestep for a low degree polynomial, we can quickly compute the higher degree optimal ones.
		\begin{table}[!ht]
			\def\arraystretch{1.5}
			\centering
			\begin{tabular}{l?{2}c|c|c}
				& \multicolumn{3}{c}{$t[s]$} \\
				Stages $S$ & $1^\text{st} \text{ Order}$ & $2^\text{nd} \text{ Order}$ & $3^\text{rd} \text{ Order}$\\
				\Xhline{5\arrayrulewidth}
				16 (Alg. from \cite{ketcheson2013optimal}) & 15.531 & 15.346 & 15.076 \\
				32 (\cref{alg:RootOpt}) & 0.478 & 0.879 & 1.183 \\
				64 (\cref{alg:RootOpt}) & 1.551 & 2.869 & 2.276 \\
				128 (\cref{alg:RootOpt}) & 9.383 & 8.845 & 8.549\\
				\Xhline{5\arrayrulewidth}
				Cumulative $t$ & 26.942 & 27.939 & 27.084
			\end{tabular}
			\caption[Runtimes of feasibility problems for different degrees constrained to spectrum of Burgers equation]{Runtimes of feasibility problems for different degrees constrained to spectrum of Burgers equation. 
				The times are obtained on a Dell Precision 5570 laptop equipped with an Intel i7-12800H.}
			\label{tab:TimesBurgers}
		\end{table}
		\item 1D Shallow Water equations: Here we optimize polynomials of degrees $30, 60$, and $120$ constrained by $\widetilde{M} = 117$ eigenvalues, again in a sequential manner.
		The pseudo-extrema and contours of the optimal $P_{60,2}(z)$ stability polynomial are displayed in \cref{fig:PE_ConvexSpectraShallowWater}.
		The CPU times are very similar to Burgers' equation with total $30.687s$ for the computation of the first order accurate stability polynomial and $28.726s$ for the second order case.
		\item 1D Ideal compressible \ac{mhd}: Polynomials of degrees $28, 56$, and $112$ are found satisfying the $\widetilde{M} = 252$ stability constraints due to the eigenvalues.
		The pseudo-extrema and contours of the optimal $P_{56,2}(z)$ stability polynomial are displayed in \cref{fig:PE_ConvexSpectraMHD}.
		The generation of the $112$ degree stability polynomial takes $94.433s$ and $98.377s$ for first and second order case, respectively.
		\item 2D Compressible Euler equations: To conclude this result section, we generate stability polynomials with degrees $26, 52$, and $104$ for a reduced set of $\widetilde{M} = 717$ eigenvalue constraints.
		The pseudo-extrema and contours of the optimal $P_{52,2}(z)$ stability polynomial are displayed in \cref{fig:PE_ConvexSpectra2D_CEE}.
		Due to the quadratic scaling in the number of constraints we have significantly longer runtimes then in the previous cases, with total $2.456 \text{ min}$ for the first order and $8.55 \text{ min}$ for the second order case.
		Note that these runtimes can be significantly reduced if only a small portion of the eigenvalues alongside the convex hull thereof is supplied as constraints.
	\end{enumerate}
	\section{Non-Convex Spectra}
	\label{sec:NonConvexSpectra}
	While a rich set of spatial discretization, equation, boundary and initial condtition combinations leads to strictly convex spectra, especially in one spatial dimension, at least equally many cases with non-convex spectra are observed.
	In this section, we propose two possible treatments thereof. 
	
	First, alpha shapes \cite{edelsbrunner1983shape, edelsbrunner2010alpha} are proposed as a potentially accurate, yet more expensive candidate which come with an unknown hyperparameter.
	Second, we follow the previous sections and use a convex hull approach also for non-convex spectra.
	\subsection{Alpha Shapes}
	\label{subsec:AlphaShapes}
\begin{figure}[!t]
		\centering
		\subfloat[$\partial \mathcal S_{50}$ for the scaled spectrum $\Delta t_{50,2} \sigma$.]{
			\label{fig:AlphaShapesStabBnd}
			\centering
			\includegraphics[width=0.45\textwidth]{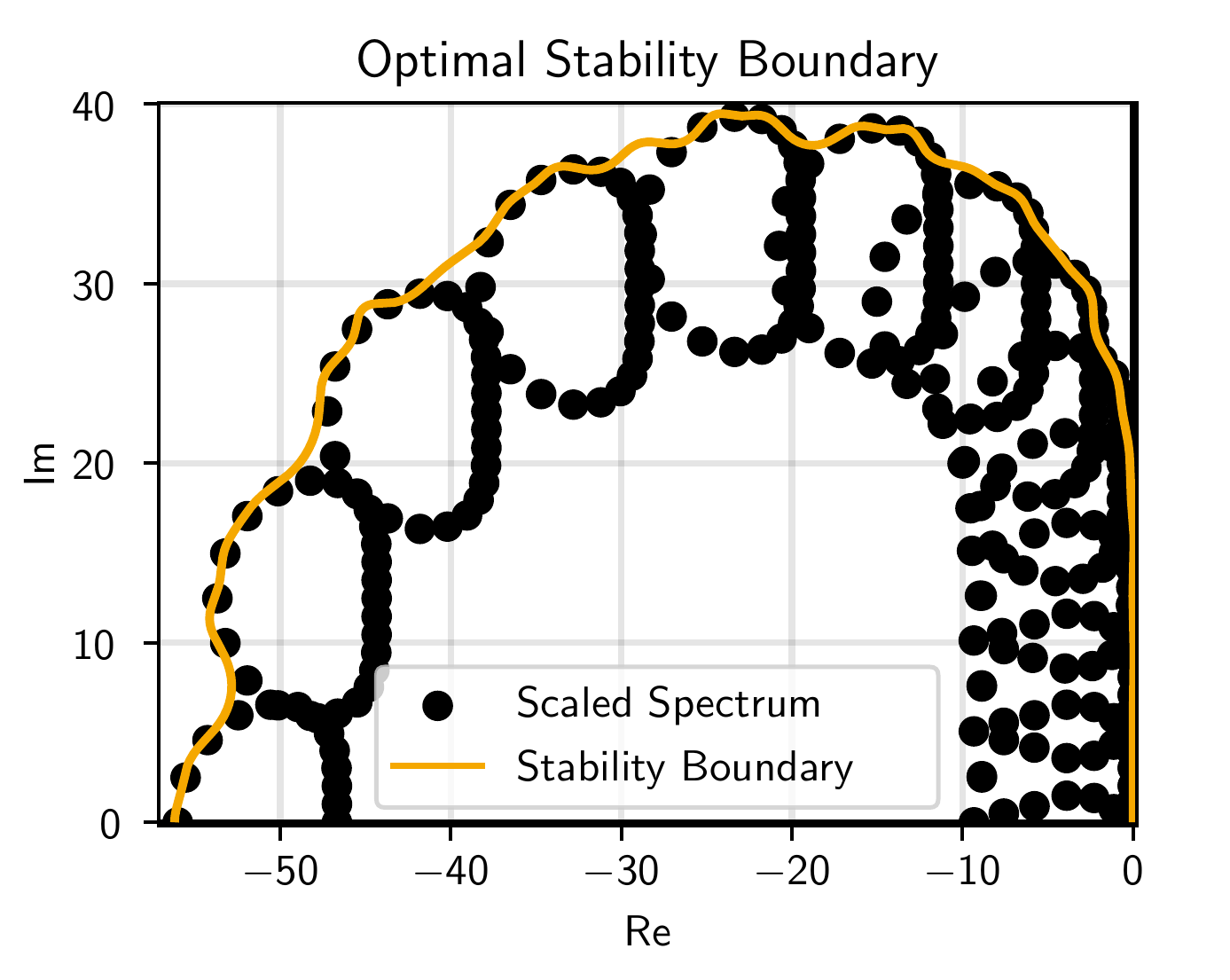}
		}
		\hfill
		\subfloat[Alpha shape with $\alpha = 0.2$ of the scaled spectrum $\Delta t_{50,2} \sigma$.]{
			\label{fig:AlphaShapes02}
			\centering
			\includegraphics[width=0.45\textwidth]{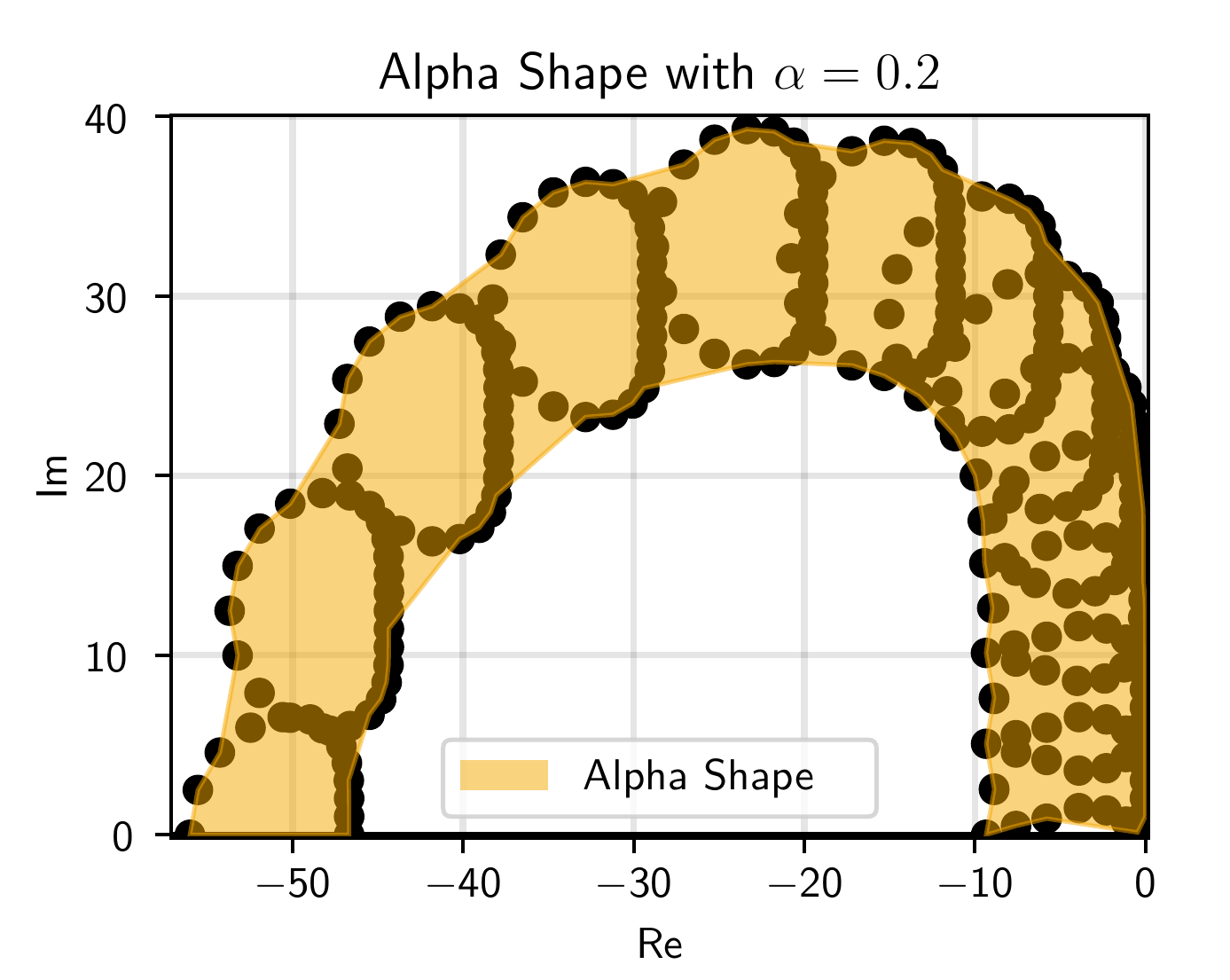}
		}
		\caption[Approximation of the Optimal Stability Boundary by the Contour of an alpha Shape]{Approximation of the optimal stability boundary by the contour of an alpha shape.
		The spectrum $\sigma$ is obtained from a \ac{dgsem} discretization of the 2D linear advection equation with velocities
		$a_x = 0.5, a_y = -0.1$ on the periodic domain $\Omega = [-1,1]^2$ discretized by a $6\times 6$ mesh with local polynomials of degree $3$ and Rusanov/Local Lax-Friedrichs flux.
	}
		\label{fig:AlphaShapes}
	\end{figure}

	Alpha shapes \cite{edelsbrunner1983shape, edelsbrunner2010alpha} provide a rigorous way of constructing a point enclosing shape that is potentially closer to the point set than the convex hull.
	Alpha shapes are parametrized by the real scalar $\alpha$ parameter that defines how close the alpha shape is shrinked to the points.
	For $\alpha > 0$, the construction of alpha shapes can intuitively described by rolling a disk with radius $1/\alpha$ around the points and drawing a line between two points of the set whenever they are on the edge of the disk and no points are in the interior of the disk.
	For $\alpha = 0$, the disk has infinite radius, i.e., degenerates to a straight line and the resulting shape is equivalent to the convex hull which is constructed by tilting the straight line around the points.
	While alpha shapes provide potentially a very accurate ansatz for the stability boundary, cf. \cref{fig:AlphaShapes} they come with two main drawbacks.
	
	First, alpha shapes introduce a hyperparameter whose choice is not clear.
	While in general higher degree polynomials follow the stability boundary more closely, i.e., correspond to larger alpha values the best choice is highly dependent on the individual spectra.
	Furthermore, for high $\alpha$ corresponding to small radii there is always the danger of obtaining a disjoint alpha shape which is not useful in this context.
	Both of these issues can be in principle addressed through manual inspection, which is clearly not attractive.
	
	Secondly, and more severely, alpha shapes provide in general only a spectrum enclosing curve $\gamma_\alpha(\tau): \mathbb R \to \mathbb C$, i.e., no spectrum enclosing function.
	This implies that for both the imaginary part (cf. \eqref{eq:PieceWiseLinearIntPolAlphaShapeImag}) and real part
	%
	\begin{align}
		\begin{split}
			\label{eq:PieceWiseLinearIntPolAlphaShape}
			I_\text{Re}\left(\tau; \sigma^\alpha_{S,p} \right) = \:& \text{Re}\left(\mu_{S,p}^{(j)} \right) + \frac{\text{Re}\left(\mu_{S,p}^{(j+1)}\right) -\text{Re}\left(\mu_{S,p}^{(j)}\right)}{\tau^{(j+1)}_{S,p} -\tau_{S,p}^{(j)}} \left(\tau - \tau_{S,p}^{(j)} \right), \\
			&\tau_{S,p}^{(j)} \leq \tau < \tau_{S,p}^{(j+1)}
		\end{split}
	\end{align}
	of the pseudo-extrema interpolation has to be performed.
	Here, $\mu^{(j)} \in \sigma^\alpha$ are the eigenvalues selected as the upper part of the alpha shape with corresponding arc length parameter $0 \leq \tau^{(j)} \leq 1$.
	This results in doubled computational costs and an overall more difficult optimization task.
	
	Nevertheless, we stress that the general idea of distributing the pseudo-extrema with equal arc length, i.e., equal distances $\Delta \tau = \tau^{(j+1)} - \tau^{(j)}, \: \forall j$ is still applicable and \cref{alg:RootOpt} can be used under the appropriate changes.
	\subsection{Convex Hull Ansatz}
	\label{subsec:ConvexHullAnsatz}
	By constructing the convex hull for non-convex spectra the framework prepared in the previous sections is recovered.
	Precisely, even if the hull of spectrum itself may be nonconvex, the convex hull is constructed and used as the ansatz for the initial placement of the pseudo-extrema.
	Using the convex hull difficulties with additional hyperparameters or the danger of disjoint alpha shapes are ruled out.
	This comes with the drawback that we can no longer expect to actually realize the maximum possible optimal timestep according to \cref{eq:TimestepLinearScaling}.
	Consequently, one has to indeed perform the optimization problems \cref{eq:OptProbx}, \cref{eq:OptProbxy} as the feasibility problems with fixed timestep will fail.
	Depending on the spectrum there are a priori no guarantees how efficient this approach will be.
	In practice, however, we find for many real-world spectra excellent results, as illustrated in the next section.
	\section{Optimal Stability Polynomials for Nonconvex Spectra}
	\label{sec:OptStabPolysNonConvex}
\begin{figure}[!t]
		\centering
		\subfloat[{Spectrum of Burgers' equation computed for the semidiscretization as specified by item \ref{item:BurgersNon}.}]{
			\label{fig:NonConvexSpectraBurgers}
			\centering
			\includegraphics[width=0.45\textwidth]{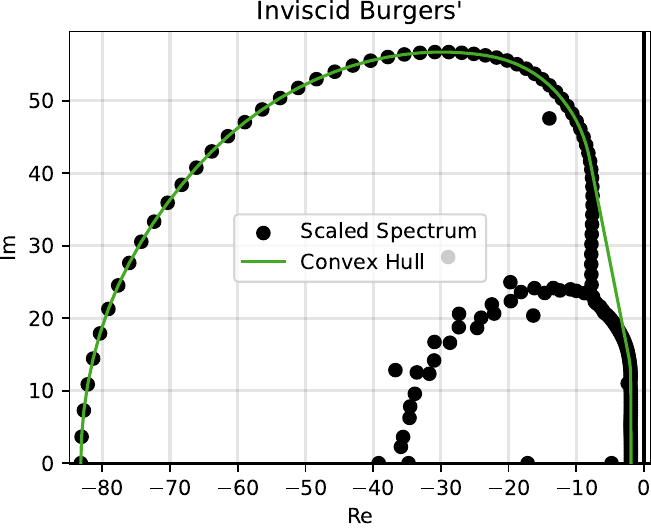}
		}
		\hfill
		\subfloat[{Spectrum of the 2D Compressible Euler equations computed for the semidiscretization as specified by item \ref{item:2D_CEE_Non}.}]{
			\label{fig:NonConvexSpectra2D_CEE}
			\centering
			\includegraphics[width=0.45\textwidth]{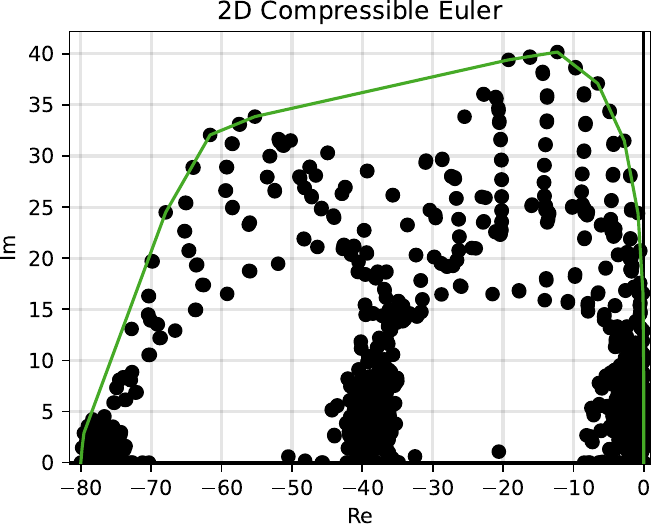}
		}
		\caption[Collection of Nonconvex Spectra]{Collection of nonconvex convex spectra of nonlinear hyperbolic \ac{pdes} and corresponding convex hull.
			The spectra are scaled with the optimal timesteps $\Delta t_{S, 2}$.
		}
		\label{fig:NonConvexSpectra}
	\end{figure}

	Here, we present results for two relevant, exemplary nonconvex spectra for which we use the convex hull-based approach.
	\begin{enumerate}
		\item \label{item:BurgersNon} Burgers' Equation with Riemannian initial data $u_0(x) = \begin{cases} 1.5 & x < 0.5 \\ 0.5 & x \geq 0 \end{cases}$.
		As before, we discretize $\Omega = [0, 1]$ with $64$ cells and use Godunov's flux.
		At $x=0$ we fix $u(t, 0) = 1.5$ while at $x=1$ an outflow boundary is used.
		The solution is reconstructed using fourth order local polynomials.
		The corresponding spectrum is displayed in \cref{fig:NonConvexSpectraBurgers} alongside its convex hull.
		The convex hull is in most parts extremely close the spectrum and we thus expect maximal timesteps very close to the optimal ones according to \cref{eq:TimestepLinearScaling}.
		Indeed, for $30, 60$, and $120$ degree polynomials we can realize the optimal timestep.
		As for the convex spectra, we display the $S=60$ case in \cref{fig:NonConvexSpectraBurgers}.
		\item \label{item:2D_CEE_Non} 2D compressible Euler equations: Isentropic vortex advection. 
		The advection of an isentropic vortex is a classic, fully nonlinear testcase with known analytical solution \cite{shu1988efficient, wang2013high}.
		In terms of the physical parameters we use the same setup as in \cite{VERMEIRE2019465, HEDAYATINASAB2022111470, VERMEIRE2021110022} which is spelled out in \cref{subsec:2DCEEIsentropicVortexAdv}.
		The numerical scheme is composed of Rusanov/Local Lax-Friedrichs flux, second order local polynomials and 8 cells in each direction.
		This little resolution is required to execute the eigenvalue decomposition in reasonable time.
		The corresponding spectrum is shown in \cref{fig:NonConvexSpectra2D_CEE} alongside its convex hull.
		Similar to Burgers' Equation, the spectrum is in most parts reasonably close to the spectrum.
		For $S=28$ we can realize about $97.9\%$ of the theoretically optimal timestep and for $S=56, S=112$ indeed the optimal timestep according to the linear timestep scaling \cref{eq:TimestepLinearScaling}.
	\end{enumerate}
	\begin{figure}[!t]
		\centering
		\subfloat[{Spectrum of Burgers' Equation computed \\for the semidiscretization defined by item \ref{item:BurgersNon}.}]{
			\label{fig:PE_NonConvexSpectraBurgers}
			\centering
			\includegraphics[width=0.45\textwidth]{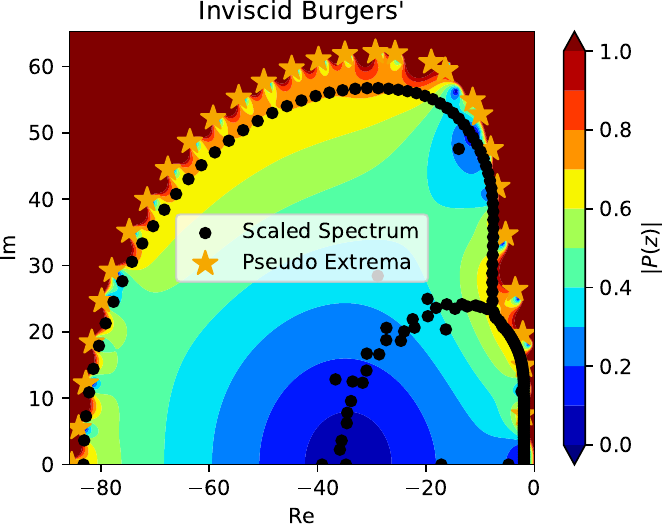}
		}
		\subfloat[{Spectrum of the 2D Compressible Euler equations computed for the semidiscretization defined by item \ref{item:2D_CEE_Non}.}]{
			\label{fig:PE_NonConvexSpectra2D_CEE}
			\centering
			\includegraphics[width=0.45\textwidth]{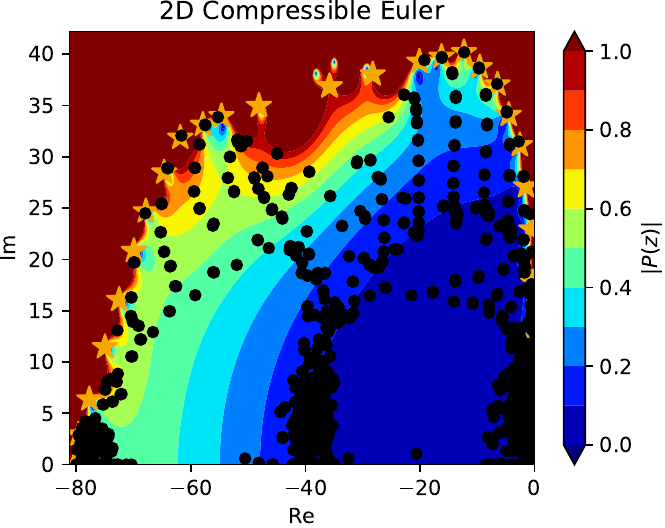}
		}
		\caption[Collection of Convex Spectra]{Collection of strictly convex spectra of nonlinear hyperbolic \ac{pdes} and optimized pseudo-extrema.
			The spectra are scaled with the optimal timestep $\Delta t_{16, 2}$ in each case.
		}
		\label{fig:PE_NonConvexSpectra}
	\end{figure}
	\section{Construction of Many Stage Runge-Kutta Methods}
	\label{sec:Construction}
	As discussed before, the primal use case of the high degree stability polynomials is the possibility to construct many-stage methods for e.g. partitioned multirate Runge-Kutta methods such as the \ac{PERK} \cite{VERMEIRE2019465, HEDAYATINASAB2022111470, VERMEIRE2023112159} schemes.
	Nevertheless, we seek to demonstrate the capabilities of the optimized stability polynomials by constructing many-stage standalone methods.
	Here, we limit ourselves to second order accurate methods for which the linear order constraints \cref{eq:ConsistencyRequirement} imply second order convergence even in the nonlinear case, without any additional constraints \cite{HairerWanner1, wanner1996solving}.
	Furthermore, we leave the construction of \ac{ssp} methods for future work and focus here on hyperbolic \ac{pdes} with smooth solutions.
	
	Due to the factorized form of the stability polynomial one can directly read-off a possible choice of the intermediate steps:
	\begin{equation}
		\label{eq:FactorizedRKMConstruction}
		P_S(z) = 1 + z \widetilde{P}_{S-1} = 
		\underbrace{1 + z \left( \prod_{j=1}^{S_\text{Real}} \underbrace{\left(1 - \frac{z}{\widetilde{r}_j} \right)}_{\text{Forward Euler}} \cdot \prod_{j=1}^{S_\text{Complex}/2} \underbrace{\left( 1- \frac{z}{\widetilde{r}_j} \right) \left( 1- \frac{z}{\widetilde{r}_j^\star} \right)}_{\text{Two-Stage Submethod}} \right)}_{\text{Forward Euler}}
	\end{equation}
	Note that we group the complex-conjugated pseudo-extrema to avoid complex-valued timesteps.
	We focus on real timesteps since they allow easy implementation into existing codes and spare us from technical complications.
	Readers interested in complex-valued timesteps are pointed to  \cite{george2021walking} and references therein.
	
	Being equipped with a stability polynomial, one has now in principle infinitely many possible choices for the actual Runge-Kutta method.
	It is well-known that internal stability, i.e., the propagation of round-off errors is a relevant concern for many stage explicit methods \cite{van1980internal, verwer1990convergence, ketcheson2014internal}.
	This issue arises also here and seems to be the major difficulty in realizing the maximum timestep.
	
	Given the factorized form of the stability polynomial it is natural to construct Runge-Kutta methods with a sequential structure.
	Here, we specify the methods in \textit{modified Shu-Osher form} \cite{ferracina2005extension, ketcheson2014internal}
	\begin{subequations}
		\label{eq:ModifiedShuOsherForm}
		\begin{align}
			\label{eq:ShuOsherUpdateEq}
			\boldsymbol Y_k & \coloneqq v_k \boldsymbol U_n + \sum_{l=1}^{k-1} \big(\alpha_{k,l} \boldsymbol Y_l + \Delta t \beta_{k,l} \boldsymbol F(\boldsymbol Y_l)\big), \quad k = 1, \dots, S+1 \\
			\boldsymbol U_{n+1} &\,= \boldsymbol Y_{S+1}
		\end{align}
	\end{subequations}
	and set $v_1 = 1, v_k = 0, k > 1$ which implies $\boldsymbol Y_1 = \boldsymbol U_n$.
	Following \cite{ketcheson2014internal} we introduce the notation 
	\begin{align}
		\alpha_{1:S}\vert_{i,j} \coloneqq \alpha_{i,j}, 
		\quad \beta_{1:S}\vert_{i,j} \coloneqq \beta_{i,j} 
	\end{align}
	with $\alpha_{1:S}, \beta_{1:S} \in \mathbb R^{S \times S}$ resembling the coefficients of the intermediate stages.
	The parameters of the final stage are denoted via $\boldsymbol \alpha_{S+1}, \boldsymbol \beta_{S+1} \in \mathbb R^{1\times S}$, i.e.,
	\begin{equation}
		\boldsymbol \alpha_{S+1} \coloneqq \begin{pmatrix}
			\alpha_{S+1,1} & \dots & \alpha_{S+1, S}
		\end{pmatrix} 
		, 
		\quad 	
		\boldsymbol \beta_{S+1} \coloneqq \begin{pmatrix}
			\beta_{S+1,1} & \dots & \beta_{S+1, S}
		\end{pmatrix} 
	\end{equation}
	which are in the present case conveniently chosen as 
	\begin{equation}
		\label{eq:alpha_beta_Splus}
		\boldsymbol \alpha_{S+1} = \begin{pmatrix}
			1 & \dots & 0
		\end{pmatrix} 
		, 
		\quad 	
		\boldsymbol \beta_{S+1} = \begin{pmatrix}
			0 & \dots & 0 & z
		\end{pmatrix}
	\end{equation}
	due to the particular representation of the stability polynomial \cref{eq:FactorizedRKMConstruction}.
	
	The two-stage submethod (cf. \cref{eq:FactorizedRKMConstruction}) corresponding to pseudo-extrema $\widetilde{r}_j, \widetilde{r}_j^\star$ are naturally represented via the low-storage scheme
	\begin{subequations}
		\label{eq:IntermediateTwoSubstage}
		\begin{align}
			\boldsymbol Y_k &= \boldsymbol Y_{k-1} + \beta_{k,k-1}^{(j)} \boldsymbol F\left(\boldsymbol Y_{k-1}\right) \\
			\boldsymbol Y_{k+1} &= \alpha_{k+1, k-1}^{(j)} \boldsymbol Y_{k-1} +  \alpha_{k+1, k}^{(j)} \boldsymbol Y_{k} + \Delta t \left( \beta_{k+1, k-1}^{(j)} \boldsymbol F(\boldsymbol U_{k-1}) + \beta_{k+1,k}^{(j)} \boldsymbol F(\boldsymbol Y_{k}) \right).
		\end{align}
	\end{subequations}
	The Forward Euler steps are readily reflected by $\alpha_{k,k-1}^{(j)} = 1, \beta_{k,k-1}^{(j)} = \frac{-1}{\widetilde{r}_j}$.
	To meet the stability polynomial \cref{eq:FactorizedRKMConstruction} we have for the two-stage methods the constraints
	\begin{subequations}
		\label{eq:ConstraintsOrder}
		\begin{align}
			\frac{ -2 \text{Re} (\widetilde{r}_j) }{\vert \widetilde{r}_j \vert^{2}} &\overset{!}{=} \alpha_{k+1,k}^{(j)} \beta_{k,k-1}^{(j)} + \beta_{k+1,k-1}^{(j)} + \beta_{k+1,k}^{(j)} \\
			\frac{1 }{\vert \widetilde{r}_j \vert^{2}} &\overset{!}{=} \beta_{k,k-1}^{(j)} \beta_{k+1,k}^{(j)}
		\end{align}
	\end{subequations}
	besides the usual requirement for consistency $\sum_{l=0}^{k-1} \alpha_{k,l} = 1$ \cite{gottlieb2001strong}.
	The timestamps/abscissae are computed via $\boldsymbol c \coloneqq (I - \alpha_{1:S})^{-1} \beta_{1:S} \boldsymbol 1$ \cite{ketcheson2014internal} where $\boldsymbol 1$ denotes a column vector of ones.
	The ordering of the two-stage submethods, i.e., the relation of $j$ to $k$ is determined based on the accumulated $\beta$ values 
	\begin{equation}
		\label{eq:Accbeta}
		\left \Vert \boldsymbol \beta^{(j)}_k \right \Vert_1 = \left\Vert \left(\beta_{k,k-1}^{(j)}, \beta_{k+1,k-1}^{(j)}, \beta_{k+1,k}^{(j)} \right) \right\Vert_1.
	\end{equation}
	The methods are then ordered such that $\left \Vert \boldsymbol \beta^{(j)}_k \right \Vert_1$ increases with $k$. 
	\subsection{Internal Stability}
	\label{subsec:InternalStabiliy}
	For the discussion of internal stability we follow the main results derived in \cite{ketcheson2014internal}.
	Under the influence of round-off errors $\boldsymbol e_k$ the modified Shu-Osher form \cref{eq:ShuOsherUpdateEq} reads
	\begin{equation}
		\label{eq:ShuOsherFormPerturbed}
		\widetilde{\boldsymbol Y}_k  \coloneqq v_k \widetilde{\boldsymbol U}_n + \sum_{l=1}^{k-1} \Big( \alpha_{j,l} \widetilde{\boldsymbol Y}_l + \Delta t \beta_{k,l} \boldsymbol F \left (\widetilde{\boldsymbol Y}_l \right) \Big) + \boldsymbol e_k, \quad k = 1, \dots, S+1 
	\end{equation}
	where we denote the perturbed, round-off error bearing stages by $\widetilde{\boldsymbol Y}_k$.
	In \cite{ketcheson2014internal} it was shown that the defect from actually computed, perturbed solution $\widetilde{\boldsymbol U}_n$ to the true solution $\boldsymbol U(t_n)$ can for linear \ac{odes} be estimated as
	\begin{equation}
		\label{eq:RoundOffErrorIneq}
		\left \Vert \widetilde{\boldsymbol U}_{n+1} - \boldsymbol U(t_{n+1}) \right \Vert \leq \left \Vert \widetilde{\boldsymbol U}_n - \boldsymbol U(t_n) \right \Vert + \sum_{k=1}^{S+1} \left \Vert Q_j(z) \right \Vert \Vert \boldsymbol e_k \Vert + \mathcal{O} \left (\Delta t^{p+1} \right).
	\end{equation}
	In addition to the usual truncation error $\mathcal{O}(\Delta t^{p+1})$ we have the amplification of round-off errors $\sum_k\Vert Q_k(z) \Vert \Vert \boldsymbol e_k \Vert$.
	Internal stability becomes a concern when both are of similar magnitude \cite{ketcheson2014internal}.
	
	In the estimate of the error of the perturbed iterates \cref{eq:RoundOffErrorIneq} $ Q_j(z)$ denote the \textit{internal stability polynomials} \cite{verwer1990convergence, ketcheson2014internal} which are for a method in Shu-Osher form computed as \cite{ketcheson2014internal}
	\begin{equation}
		\boldsymbol Q(z; \alpha, \beta)
		\coloneqq (\boldsymbol \alpha_{S+1} + z \boldsymbol\beta_{S+1} ) \left( I - \alpha_{1:S} - z \beta_{1:S} \right)^{-1}.
	\end{equation}
	For explicit methods 
	there is no initial round-off error $\boldsymbol e_1$ \cite{ketcheson2014internal} and it suffices to consider the internal stability polynomials starting from second stage.
	To estimate the potential of round-off error amplification it is customary to investigate 
	\begin{equation}
		\label{eq:MaxAmpFactorAlternative}
		\widetilde{\mathcal{M}}(\alpha, \beta) \coloneqq \max_{z \in \mathcal S}\sum_{k=2}^{S+1} \vert Q_k(z; \alpha, \beta) \vert
	\end{equation}
	which is a more precise variant of the \textit{maximum internal amplification factor} proposed in \cite{ketcheson2014internal}.
	For standard double precision floating point datatypes we expect the round-off errors to be $\varepsilon = \mathcal{O}\left( 10^{-15}\right)$ at most, thus we expect the overall internal errors to be of order $\widetilde{\mathcal{M}} \cdot 10^{-15}$.
	This can then be compared to $\mathcal{O}(\Delta t^{p+1})$ yielding an approximate criterion to estimate the influence of round-off errors.
	Both metrics are reported for the examples considered in \cref{sec:ResultsStrictlyConvexSpectra}.
	
	In principle, one could now set out to minimize $\widetilde{\mathcal M}$ over $\alpha, \beta$ under the constraints \cref{eq:ConstraintsOrder}.
	Given that even for methods with two intermediate stages \cref{eq:IntermediateTwoSubstage} we have $\frac{5}{2} S_\text{Complex} + S_\text{Real}$ free optimization variables and one is minimizing not only one, but $S$ polynomials at $M$ eigenvalues this is a significantly harder optimization problem than the original one \cref{eq:OptProb}.
	In addition, there are no immediate properties of the internal stability polynomials $Q_j$ available, rendering this problem practically infeasible.
	As an alternative, we propose a significantly cheaper heuristically motivated approach in the next section.
	\subsection{Internal Stability Optimization: Heuristic Approach}
	We begin by recalling that $\alpha_{1:S}, \beta_{1:S}$ are for explicit methods strictly lower triangular, i.e., nilpotent matrices with index $S$.
	Consequently, the identity 
	\begin{equation}
		\left( I - \alpha_{1:S} - z \beta_{1:S} \right)^{-1} = \sum_{k=0}^{S-1} \left(\alpha_{1:S} + z \beta_{1:S}\right)^k
	\end{equation}
	holds and allows an alternative to the inversion of the typically ill-conditioned matrix.
	We make now the following heuristic argument:
	For $\alpha_{k,l} \in [0,1]$ and $\vert z \vert \sim \mathcal{O}(100)$ for the many stage methods we expect a significantly higher sensitivity of $\widetilde{M}$ with respect to $\beta_{k,l}$ than $\alpha_{k,l}$.
	In order to minimize $\widetilde{M}$ we thus minimize the accumulated $\beta_{k,l}$ values \cref{eq:Accbeta} for each two-stage submethod individually.
	For these optimization problems constraints \cref{eq:ConstraintsOrder} satisfying local minima are easily found.
	Note that negative $\beta_{k,l}$ are allowed here, which allow further reduction of $\widetilde{M}$ due to the increased flexibility.
	To illustrate this, consider the $104$ degree stability polynomial optimized for the scalar advection equation (see \cref{subsec:1D_Adv_Gauss}) for which we construct Runge-Kutta methods with non-negative and unconstrained $\beta_{k,l}$.
	For this case, values of $\mathcal M$ for grouped and non grouped pseudo-extrema (see next section) with non-negative and unconstrained $\beta_{k,l}$ are tabulated in \cref{tab:MNegBetaGrouping}.
	\subsection{Generalized Lebedev's Idea}
	For the higher order \ac{dg} spectra we observe a couple of pseudo-extrema near the imaginary axis with $\left \vert \text{Im}(\widetilde{\boldsymbol r}_j) \right \vert \gg -\text{Re}(\widetilde{\boldsymbol r}_j)$ which can lead to $\beta_{k,l} \gg 1$, which should be avoided.
	Similar to Lebedev's Idea/Realization \cite{lebedev1989explicit, lebedev2017solve} we group the complex conjugated pseudo-extrema with $\text{Re}(\widetilde{\boldsymbol r}_j) > -0.5$ with the pseudo-extrema of largest $\left \vert \text{Re}(\widetilde{\boldsymbol r}_j) \right \vert$.
	By staying with the two register form \cref{eq:IntermediateTwoSubstage} we have now $13$ free parameters, three linear constraints due to $\sum_{k=0}^{j-1} \alpha_{j,k} = 1$ and four nonlinear constraints which are derived similar to \cref{eq:ConstraintsOrder} by multplying out polynomials.
	As they are quite lengthy but readily obtained we omit them here.
	As for the remaining two-stage methods, we optimize $\alpha, \beta$ to minimize $\left \Vert \boldsymbol \beta^{(j)}_k \right \Vert_1$.
	
	The grouping of pseudo-extrema with small and large real part has a dramatic influence on the internal stability properties.
	For instance, consider again the $104$ degree, second order accurate stability polynomial optimized for the scalar advection equation (see \cref{subsec:1D_Adv_Gauss}) for which we construct Runge-Kutta methods with and without grouping.
	Values for $\widetilde{\mathcal M}$ are tabulated in \cref{tab:MNegBetaGrouping} highlighting the importance of the grouping.
	Due to this results all in the following constructed methods bear grouping of pseudo-extrema and may have negative $\beta_{k,l}$.
	\begin{table}[!ht]
		\def\arraystretch{1.5}
		\centering
		\begin{tabular}{l?{2}l|l}
			& $\beta_{k,l} \in \mathbb R_0^+$ & $\beta_{k,l} \in \mathbb R$ \\
			\Xhline{5\arrayrulewidth}
			No Grouping & $2.39 \cdot 10^{22}$ & $2.19 \cdot 10^{22}$ \\
			Grouping & $7.54 \cdot 10^{16}$ & $1.63 \cdot 10^{9}$ \\
		\end{tabular}
		\caption[Values of Internal Error Amplification for different Method Construction Paradigms]{Internal error amplification $\widetilde M$ \cref{eq:MaxAmpFactorAlternative} for different method construction paradigms.}
		\label{tab:MNegBetaGrouping}
	\end{table}
	\subsection{Note on \ac{ssp} Properties}
	It is natural to ask for \ac{ssp} properties of Runge-Kutta schemes intended for the integration of (nonlinear) hyperbolic systems.
	Furthermore, given the representation of the method in Shu-Osher form \cref{eq:ModifiedShuOsherForm} one can readily tell whether the method is monotonicity preserving or not.
	For the methods resembling the factorized form of the stability polynomial the coefficients of the last stage are chosen according to \cref{eq:alpha_beta_Splus}.
	Consequently, the \ac{ssp} cofficient \cite{gottlieb2001strong} 
	\begin{equation}
		c \coloneqq \min_{k,l} \frac{\alpha_{k,l}}{\beta_{k,l}}
	\end{equation}
	is due to $\alpha_{S+1,S} =0, \beta_{S+1,S} = 1$ always zero, regardless of the parametrization of the intermediate stages.
	This implies that there exists no timestep $\Delta t >0$ such that the Runge-Kutta method is guaranteed to be monotonicity preserving.
	One will also need to revisit whether negative $\beta_{k,l}$ should be allowed as they demand special treatment for \ac{ssp} methods \cite{gottlieb2001strong}.
	The construction of many stage \ac{ssp} methods is ongoing research. 
	\section{Many-Stage Runge-Kutta Methods: Results}
	\label{sec:ResultsStrictlyConvexSpectra}
	We present convergence studies of the many-stage Runge-Kutta methods which are constructed from the factorized form of the stability polynomial using the generalization of Lebedev's Idea and possibly negative $\beta_{k,l}$.
	\subsection{Linear Problems}
	\subsubsection{1D Advection of a Gaussian Pulse}
	\label{subsec:1D_Adv_Gauss}
	As a first example, we consider the linear advection equation with unit transport velocity and Gaussian initial data
	\begin{equation}
		\label{eq:AdvEq}
		u_t + u_x = 0, \quad u(x, 0) = u_0(x) = \exp\left( \frac{-x^2}{0.1}\right)
	\end{equation}	
	on $\Omega$ = $[-5, 5]$ equipped with periodic boundaries.
	We discretize $\Omega$ with $512$ cells and \cref{eq:AdvEq} via the \ac{dgsem} with third order polynomials and Rusanov/Local Lax-Friedrichs flux implemented in \texttt{Trixi.jl} \cite{trixi1, trixi2, trixi3}.
	Using \texttt{Trixi.jl} we can generate the Jacobian using algorithmic differentiation and compute the corresponding spectrum.
	For this spectrum we compute the reference timestep $\Delta_{16,3} \approx 3.53 \cdot 10^{-2}$ using the algorithm proposed in \cite{ketcheson2013optimal} which is then supplied to the feasibility problems \cref{eq:FeasProbx} and \cref{eq:FeasProbxy}.
	
\begin{figure}[!t]
		\centering
		\subfloat[{$L^\infty$ errors for scalar advection equation.}]{
			\label{fig:ConvAdvLInf}
			\centering
			\includegraphics[width=0.45\textwidth]{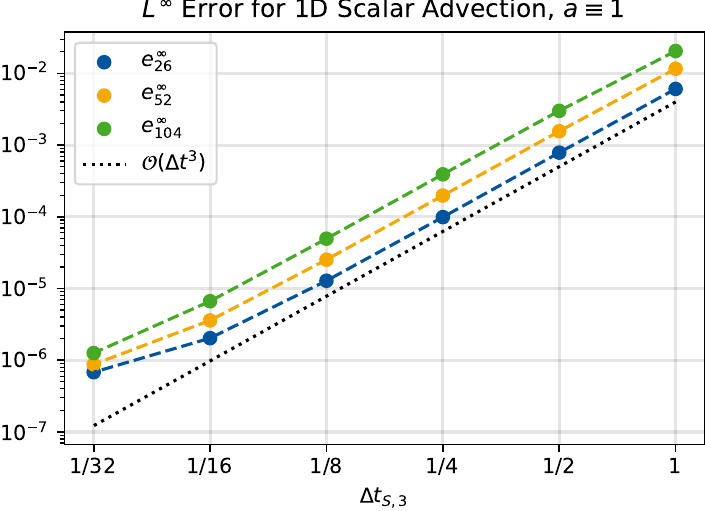}
		}
		\subfloat[{Weighted $L^1$ errors for scalar advection equation.}]{
			\label{fig:ConvAdvL1}
			\centering
			\includegraphics[width=0.45\textwidth]{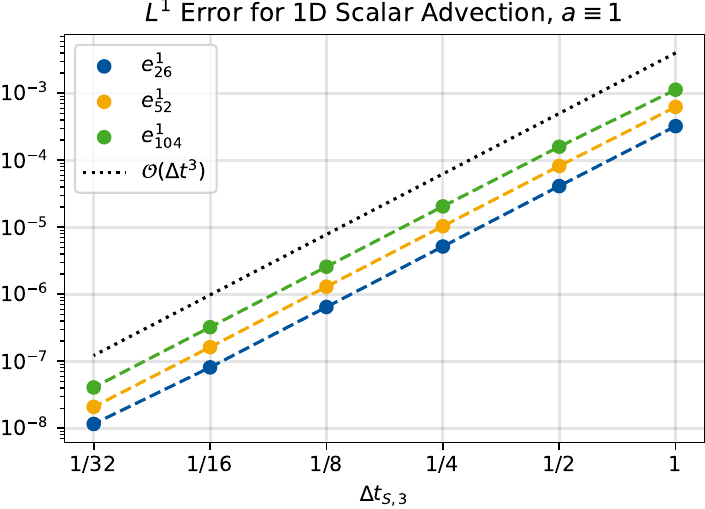}
		}
		\caption[Convergence plots for Methods optimized for 1D Scalar Advection Equation]{Third order convergence in $L^\infty, L^1$ for methods with 26, 52, and 104 stages optimized for the 1D scalar advection equation \eqref{eq:AdvEq}.
		}
		\label{fig:ConvAdv}
	\end{figure}

	We compute ten passes through the domain corresponding to final time $t_f = 100$ with $26, 52,$ and $104$ stage, third order accurate methods using the Runge-Kutta parametrization with negative $\beta$ and generalized Lebedev's idea.
	The used timestep $\Delta t$, number of taken timesteps $N_t$, internal amplification factors $\widetilde{\mathcal M}$ and the order of the truncation error $\Delta t^4$ are reported in \cref{tab:AdvGaussian}.
	The third order convergence in $L^\infty$-norm
	\begin{equation}
		e_S^\infty \coloneqq \left \Vert u^{(h,S)}(t_f, x) - u(t_f, x) \right \Vert_\infty
	\end{equation}
	is observed until spatial discretization errors are becoming relevant, cf. \cref{fig:ConvAdvLInf}.
	In weighted $L^1$-norm 
	\begin{equation}
		\label{eq:L1Norm}
		e_S^1 \coloneqq \frac{1}{\vert \Omega \vert } \left \Vert u^{(h,S)}(t_f, x) - u(t_f, x) \right \Vert_1
	\end{equation}
	third order convergence is observed except for the smallest timestep, as displayed in \cref{fig:ConvAdvL1}.
	We observe that the error constant increases with stage number, which seems to be a general phenomenon for the herein constructed methods.
	%
	%
	\begin{table}
		\def\arraystretch{1.2}
		\centering
		\begin{tabular}{c?{2}c|c|c|c}
			$S$ & $\Delta t$ & $N_t$  & $ \widetilde{\mathcal{M}} \cdot 10^{-15}$ & $\Delta t^4$ \\
			\Xhline{5\arrayrulewidth}
			26 & $5.72 \cdot 10^{-2}$ & $1748$ & $ 3.59 \cdot 10^{-12}$ & $1.07 \cdot 10^{-5}$ \\
			52 & $1.14 \cdot 10^{-1} $ & $874$ & $1.46 \cdot 10^{-9}$ & $1.72 \cdot 10^{-4}$\\
			104 & $2.29 \cdot 10^{-1} $ & $437$ & $5.39 \cdot 10^{-5}$ & $2.75 \cdot 10^{-3}$
		\end{tabular}
		\caption[Experimental data Advection Gaussian Pulse]
		{Results for the simulation of the 1D scalar advection equation \cref{eq:AdvEq} with three many stage optimized methods constructed using \cref{alg:RootOpt}.
			The timestep $\Delta t$ is kept constant over the $N_t$ timesteps, with the possible exemption of the last timestep which may be reduced in order to match the desired final time $t_f = 100$.}
		\label{tab:AdvGaussian}
	\end{table}
	\subsubsection{2D Linearized Euler Equations}
	The linearized Euler equations with source term read in two spatial dimensions, primitive variables
	\begin{equation}
		\label{eq:2DLinEuler}
		\partial_t
		\begin{pmatrix}
			\rho' \\ u' \\ v' \\ p'
		\end{pmatrix}
		+
		\partial_x
		\begin{pmatrix}
			\bar{\rho} u' + \bar{u} \rho ' \\ \bar{u} u' + \frac{p'}{\bar{\rho}} \\ \bar{u} v' \\ \bar{u} p' + c^2 \bar{\rho} u'
		\end{pmatrix}
		+
		\partial_y
		\begin{pmatrix}
			\bar{\rho} v' + \bar{v} \rho ' \\ \bar{v} u' \\ \bar{v} v' + \frac{p'}{\bar{\rho}} \\ \bar{v} p' + c^2 \bar{\rho} v'
		\end{pmatrix}
		=
		\boldsymbol s
	\end{equation}
	where $(\bar{\rho}, \bar{u}, \bar{v}, c) = (1, 1, 1, 1)$ denote the base state.
	By the method of \textit{manufactured solutions} \cite{salari2000code, roache2002code} we set the solution to
	\begin{equation}
		\label{eq:2DLinEuler_IC}
		\begin{pmatrix}
			\rho'\\ u'\\ v' \\ p'
		\end{pmatrix}
		= \begin{pmatrix}
			-\cos(2\pi t) \left[\sin(2\pi x) - \sin(2\pi y) \right] \\ \sin(2\pi t)\cos(2 \pi x) \\ \sin(2\pi t)\cos(2 \pi y) \\ -\cos(2\pi t) \left[\sin(2\pi x) - \sin(2\pi y) \right]
		\end{pmatrix}
	\end{equation}
	with corresponding source terms
	\begin{equation}
		\boldsymbol s(t,x,y) = 2\pi \begin{pmatrix}
			-\cos(2\pi t) \left[\cos(2\pi x) - \cos(2\pi y) \right] \\ \sin(2\pi t)\sin(2 \pi x) \\ \sin(2\pi t)\sin(2 \pi y) \\ -\cos(2\pi t) \left[\cos(2\pi x) - \cos(2\pi y) \right]
		\end{pmatrix}.
	\end{equation}
	The linearized Euler equations \cref{eq:2DLinEuler} are discretized on $\Omega = [0, 1]^2$ with $64$ elements in each direction.
	We employ the \ac{dgsem} with fourth order local polynomials and \ac{hll} \cite{harten1983upstream} flux on a periodic domain.
	To leverage the computational costs of the eigenvalue decomposition, we compute the spectrum for a discretization with $8$ cells in each coordinate direction.
	For the resulting set of eigenvalues we compute the reference timestep $\Delta t_{16,2} \approx 4.53 \cdot 10^{-2}$ using the algorithm proposed in \cite{ketcheson2013optimal} which is then supplied to the feasibility problems \cref{eq:FeasProbx} and \cref{eq:FeasProbxy}.
	
\begin{figure}[!t]
		\centering
		\includegraphics[width=0.45\textwidth]{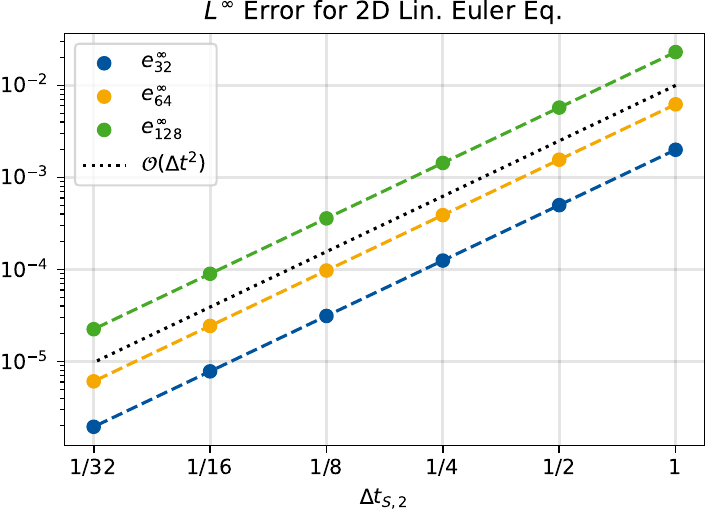}
		\caption[Second Order Convergence for the 2D Linearized Euler Equations]{Second order convergence of the optimized methods with $32, 64$, and $128$ stages optimized for the 2D linearized Euler equations.
			For each method, we reduce from the largest stable timestep $\text{CFL} \cdot \Delta t_{S,2}$.
		}
		\label{fig:Conv2D_LinEuler}
	\end{figure}

	For the actual simulation on the $8$ times finer mesh we reduce the obtained timesteps accordingly.
	A convergence study is carried out in terms of the error in density fluctuation $\rho$ at final time $t_f = 10.5$ and displayed in \cref{fig:Conv2D_LinEuler}.
	We use $t_f= 10.5$ since we observed for $t_f = 10$ spurious convergence of higher than second order, which is a peculiarity of the employed testcase rather than a general phenomenon.
	Again, second order convergence is observed in $L^\infty$ for all timesteps.
	Here, we construct polynomials of degrees $32, 64$, and $128$ and the corresponding methods with parametrization using potentially negative $\beta$ and the generalized Lebedev's idea.
	Since the right-hand-side of the \ac{ode} system \eqref{eq:ODESys2} depends also explicitly on time this example also showcases convergence of the many-stage methods for non-homogeneous problems.
	
	In contrast to the previous example, we observe that the theoretically possible maximum stable timestep has to be reduced by factor $\text{CFL}$ (see \cref{tab:2DLinEuler}) for a stable simulation which we attribute to issues with internal stability.
	Consider for instance the $128$ stage method with $\text{CFL} = 1.0$ (corresponding to an unstable simulation) with $\widetilde{\mathcal M} \cdot 10^{-15} \approx 7.60 \cdot 10^{-2}$ while $\Delta t^3 \approx 9.29 \cdot 10^{-4}$.
	As discussed in \cite{ketcheson2014internal}, we observe in this case stability problems when error terms due to round off are of  larger magnitude than the truncation error.
	For the $S=64$ stability polynomial the reduction in optimal timestep is a consequence of the optimization of the stability polynomials for the reduced spectrum, rather than related to internal stability.
	\begin{table}
		\def\arraystretch{1.2}
		\centering
		\begin{tabular}{c?{2}c|c|c|c|c}
			$S$ & $\Delta t$ & CFL & $N_t$  & $ \widetilde{\mathcal{M}} \cdot 10^{-15}$ & $\Delta t^3$ \\
			\Xhline{5\arrayrulewidth}
			32 & $1.13 \cdot 10^{-2}$ & $1.0$ & $928$ & $ 1.91 \cdot 10^{-11}$ & $1.45 \cdot 10^{-6}$ \\
			64 & $2.17 \cdot 10^{-2} $ & $0.96$ & $484$ & $6.30 \cdot 10^{-9}$ & $1.03 \cdot 10^{-5}$\\
			128 & $3.67 \cdot 10^{-2} $ & $0.81$ & $287$ & $5.63 \cdot 10^{-5}$ & $4.94 \cdot 10^{-5}$
		\end{tabular}
		\caption[Experimental data 2D linearized Euler equations]
		{Results for the simulation of the 2D linearized Euler equations \cref{eq:2DLinEuler} with three many stage optimized methods constructed using \cref{alg:RootOpt}.
			The timestep $\Delta t$ is kept constant over the $N_t$ timesteps, with the possible exemption of the last timestep which may be reduced in order to match the desired final time $t_f = 10.5$.}
		\label{tab:2DLinEuler}
	\end{table}
	\subsection{Nonlinear Problems}
	\subsubsection{Burgers' Equation with Source Term}
	We consider Burgers' Equation as the classic prototype equation for nonlinear hyperbolic \ac{pdes}.
	With source term $s(t, x)$, Burgers' Equation reads
	\begin{equation}
		\label{eq:BurgersEqSourceTerm}
		\partial_t u + \frac12 \partial_x u^2	= s.
	\end{equation}
	The presence of the source term allows a convenient construction of a continuously differentiable solution $u(t, x) \in C^1$.
	This approach is known as the method of \textit{manufactured solutions} \cite{salari2000code, roache2002code}.
	Here, the solution is set to 
	\begin{equation}
		u(t, x) \coloneqq 2 + \sin\big(2 \pi (x-t)\big)
	\end{equation}
	which corresponds to a simple advection of the initial condition, although this time through a nonlinear \ac{pde}.
	Note that $u(t, x)$ is periodic on $[0, 1]$ for all times $t$.
	The source term is then computed from the governing \ac{pde} \cref{eq:BurgersEqSourceTerm} as
	\begin{equation}
		\label{eq:BurgersSourceTermManufacturedSol}
		s(t, x) \coloneqq 2\pi \cos\big(2 \pi (x-t)\big) \Big(1 + \sin\big(2 \pi (x-t)\big) \Big).
	\end{equation}
	The spatial discretization is realized using the \ac{dgsem} with third order polynomials on a $256$ element mesh $\Omega^{(h)}$ with periodic boundaries and Godunov flux \cite{OsherRiemann1984}.
	
	Despite being effectively also only the simple advection of the initial data $u_0(x)$, it is now observed that internal stability limits the effectively stable timestep much more than in case of the advection equation discussed previously.
	This is attributed to the fact that the perturbations $\boldsymbol e_j$ do not only corrupt the stage updates (cf. \cref{eq:ShuOsherFormPerturbed}), but thereby also change the point of evaluation of the Jacobian and thus the spectrum.
	This can have dramatic effects since the stability polynomial is only optimized for the spectrum corresponding to the unperturbed solution.
	This is further amplified by the lack of the \ac{ssp} property of the constructed schemes.
	Thus, as the time integration schemes used here have no oscillation-surpressing guarantees, these do arise and further perturb the spectrum.
	To measure this, we compute the increase of the \textit{total variation}
	\begin{equation}
		e_\text{TV} \coloneqq \Vert \boldsymbol u^{(h)}\left(t_f\right) \Vert_\text{TV} - \Vert \boldsymbol u^{(h)}\left(t_0\right) \Vert_\text{TV} 
	\end{equation}
	where the total variation semi-norm is defined as usual
	\begin{equation}
		\Vert \boldsymbol u(t) \Vert_\text{TV} \coloneqq \sum_{j} \left \vert u_{j+1}(t) - u_{j}(t) \right \vert .
	\end{equation}
	For methods with $S= 28, 56$, and $112$ stages we observe indeed increase in total variation, cf. \cref{tab:BurgersSmooth}.
	\begin{table}
		\def\arraystretch{1.2}
		\centering
		\begin{tabular}{c?{2}c|c|c|c|c|c}
			$S$ & $\Delta t$ & CFL & $N_t$  & $ \widetilde{\mathcal{M}} \cdot 10^{-15}$ & $\Delta t^3$ & $e_\text{TV}$ \\
			\Xhline{5\arrayrulewidth}
			28 & $3.67 \cdot 10^{-3}$ & $0.72$ & $1363$ & $ 1.41 \cdot 10^{-12}$ & $4.95 \cdot 10^{-8}$ & $1.30 \cdot 10^{-2}$\\
			56 & $7.04 \cdot 10^{-3} $ & $0.69$ & $711$ & $1.39 \cdot 10^{-10}$ & $3.48 \cdot 10^{-7}$ & $5.18 \cdot 10^{-2}$ \\
			112 & $3.99 \cdot 10^{-3} $ & $0.28$ & $876$ & $4.79 \cdot 10^{-10}$ & $1.86 \cdot 10^{-7}$ & $6.09 \cdot 10^{-3}$
		\end{tabular}
		\caption[Experimental data Burgers' equation with source term]
		{Results for the simulation of Burgers equation \cref{eq:BurgersEqSourceTerm} with three many stage optimized methods constructed using \cref{alg:RootOpt}.
			The timestep $\Delta t$ is kept constant over the $N_t$ timesteps, with the possible exemption of the last timestep which may be reduced in order to match the desired final time $t_f = 5$.}
		\label{tab:BurgersSmooth}
	\end{table}
	Here, the significant decrease of the optimal timestep even for methods with moderate stages is not expected to follow from issues with internal stability.
	Although the estimate for the perturbed approximation \cref{eq:RoundOffErrorIneq} is only valid for linear systems, we consider it also here as a crude estimate for the influence of internal stability.
	As note in table \cref{tab:BurgersSmooth} the maximum possible timestep $\Delta t_{28,2}$ has to be decreased to $72\%$ of its theoretically possible value.
	As $\widetilde{M} \cdot 10^{-15}$ is still 4 orders of magnitude smaller than $\Delta t^3$ we do not expect the round off errors to limit the stability.
	
	Instead, we believe that the lack of the \ac{ssp} property causes the limitations in timestep.
	To provide reason for this, we consider \ac{ssp} the third order accurate methods with variable number of stages $S_n = n^2$ proposed in \cite{Ketcheson2008highly}.
	For the $S=15^2 = 225 $ method we have $\widetilde{M} \cdot 10^{-15} \approx 2.13 \cdot 10^{-13}$ and $\Delta t^4 \approx 5.42 \cdot 10^{-8}$, i.e., similar values to the $28$ stage method constructed here.
	In that case, however, we have $e_\text{TV} \equiv 0$, as expected for the \ac{tvd} time integration when applied to a smooth solution.
	For the schemes constructed here the spurious oscillations can even produce amplifying modes in the spectrum of the semidiscretization, i.e., move the system into an unphysical, energy-generating state.
	In that case, the true solution itself is unstable - this is a qualitatively different scenario to the lack of A-stability where the numerical scheme is unstable for a certain eigenvalue.
	For $S=28$ and $\Delta t$ from \cref{tab:BurgersSmooth} we have for instance at $t=0.5$ an eigenvalue with $\text{Re} \left( \lambda \right) = 8.57 \cdot 10^{-12}$ which bears the possibility to cause the whole computation to diverge.
	
\begin{figure}[!t]
		\centering
		\subfloat[{$L^\infty$ errors for Burgers' equation.}]{
			\label{fig:ConvBurgersSmoothLInf}
			\centering
			\includegraphics[width=0.45\textwidth]{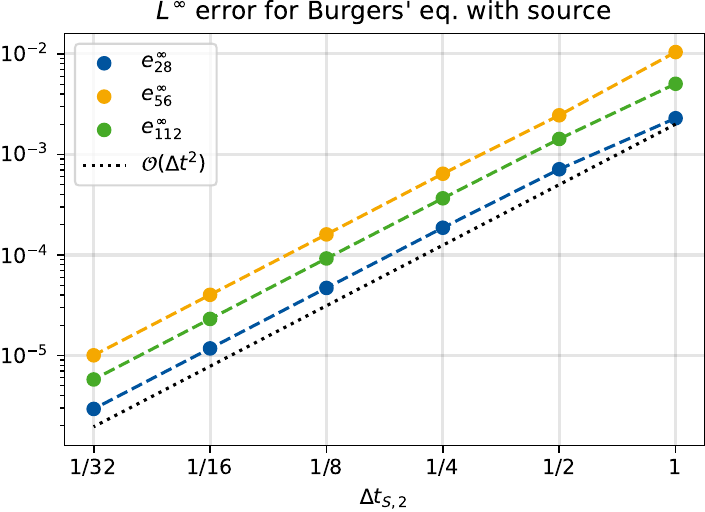}
		}
		\subfloat[{Weighted $L^1$ errors for Burgers' equation.}]{
			\label{fig:ConvBurgersSmoothL1}
			\centering
			\includegraphics[width=0.45\textwidth]{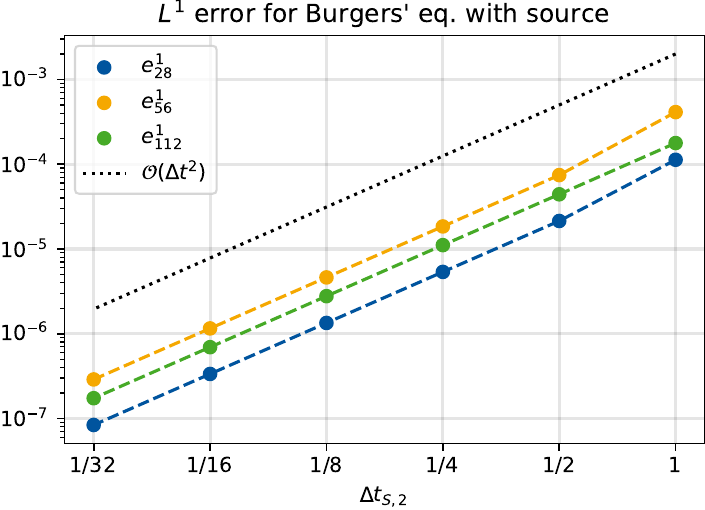}
		}
		\caption[Convergence plots for Methods optimized for Burgers' Equation with source]{Second order convergence in $L^\infty, L^1$ for methods optimized for Burgers' Equation with source \cref{eq:BurgersEqSourceTerm}.
		}
		\label{fig:ConvBurgersSmooth}
	\end{figure}

	Nevertheless, all methods achieve second order convergence in $L^\infty$-norm once the oscillations vanished which is the case for $\Delta t \leq \frac{1}{4} \Delta t_{S,2}$, cf. \cref{fig:ConvBurgersSmoothLInf}.
	In weighted $L^1$-norm \eqref{eq:L1Norm} second order convergence is observed for every timestep as displayed in \cref{fig:ConvBurgersSmoothL1}.
	\subsubsection{2D Compressible Euler Equations: Isentropic Vortex Advection}
	\label{subsec:2DCEEIsentropicVortexAdv}
	To conclude the examples, we consider the 2D compressible Euler equations 
	\begin{subequations}
		\label{eq:EulerEqs}
		\begin{align}
			\partial_t \rho + \nabla \cdot \begin{pmatrix} \rho v_x \\ \rho v_y \end{pmatrix} &= 0\\
			\label{eq:ConsMomentum}
			\partial_t \begin{pmatrix} \rho v_x \\ \rho v_y \end{pmatrix} + \nabla \cdot \begin{pmatrix}
				\rho v_x^2 & \rho v_x v_y \\ \rho v_y v_x & \rho v_y^2 
			\end{pmatrix} + \nabla p &= \boldsymbol 0 \\
			\partial_t E + \nabla \cdot \begin{pmatrix} (E + p) v_x \\ (E + p) v_y\end{pmatrix}  &=0
		\end{align}
	\end{subequations}
	with total energy $E = E(\rho, \boldsymbol v, p) = \rho \left(\frac{p}{\gamma - 1} + \frac{1}{2} \left(v_x^2 + v_y^2\right) \right)$. 
	As a testcase we consider the advection of an isentropic vortex \cite{shu1998essentially, wang2013high}.
	Here, we use similar parameters to the ones used in \cite{HEDAYATINASAB2022111470,VERMEIRE2019465}.
	In particular, the base state is set to 
	\begin{equation}
		\label{eq:IsentropicVortexAdvectionBaseState}
		\rho_\infty = 1, \quad \boldsymbol v_\infty = \begin{pmatrix}
			1 \\ 1
		\end{pmatrix}, \quad p_\infty \coloneqq \frac{\rho_\infty^\gamma}{\gamma \text{Ma}_\infty^2}
	\end{equation}
	with $\text{Ma}_\infty = 0.4$.
	To localize the effect of the vortex centered at $	\boldsymbol c(t, x, y) \coloneqq \begin{pmatrix}
		x \\ y
	\end{pmatrix} - \boldsymbol v_\infty t$
	the perturbations are weighted with the Gaussian $g(t, x, y) \coloneqq \exp\left(\frac{1 - \Vert \boldsymbol c(t, x, y) \Vert_2^2}{2R^2}\right)$
	where $R= 1.5$.
	While the size of the vortex is governed by $R$, its intensity/strength is quantified by $I$ which is set here to $13.5$ following \cite{HEDAYATINASAB2022111470,VERMEIRE2019465}.
	The density is given by 
	\begin{equation}
		\label{eq:IsentropicVortexDensity}
		\rho(t, x, y) = \rho_\infty \left( 1 - \frac{I^2 M^2 (\gamma-1) g^2(t, x, y)}{8\pi^2} \right)^\frac{1}{\gamma -1}
	\end{equation}
	and the corresponding perturbed velocities are 
	\begin{equation}
		\boldsymbol v(t, x, y) = \boldsymbol v_\infty + \frac{I g(t, x,y)}{2 \pi R} \boldsymbol c(t, x, y)
	\end{equation}
	while the pressure is computed analogous to the base pressure \cref{eq:IsentropicVortexAdvectionBaseState} as
	\begin{equation}
		\label{eq:IsentropicVortexPressure}
		p(t, x,y) =  \frac{\rho^\gamma(t, x, y)}{\gamma \text{Ma}_\infty^2}.
	\end{equation}

\begin{figure}[!t]
		\centering
		\includegraphics[width=0.45\textwidth]{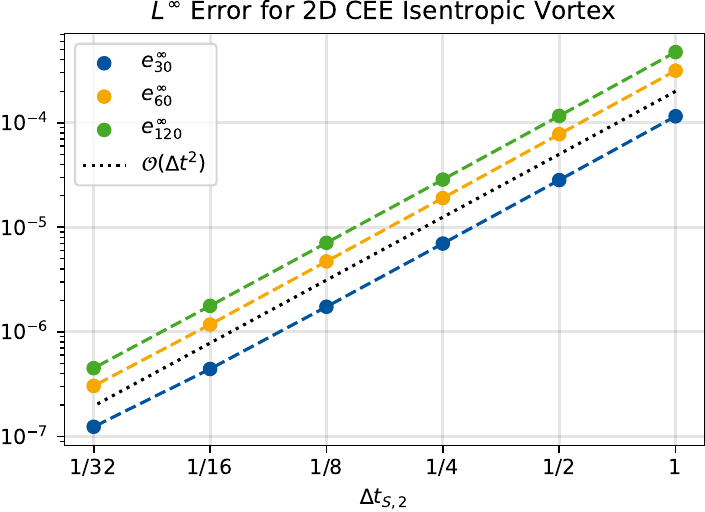}
		\caption[Second Order Convergence for the 2D Compressible Euler Equations]{Second order convergence of the optimized methods with $30, 60$, and $120$ stages optimized for the 2D compressible Euler equations.
			For each method, we reduce from the largest stable timestep $\text{CFL} \cdot \Delta t_{S,2}$.
		}
		\label{fig:2DCEE}
	\end{figure}

	The advection of the vortex is simulated on $\Omega = [-10,10]^2$ discretized with $64\times64$ elements, sixth order local polynomials, and HLLC Flux \cite{toro1994restoration}.
	Final time $t_f$ is set to $20$ which corresponds to a full traversion of the vortex through the periodic domain.
	For this setup $30,60$, and $120$ degree methods of second order accuracy are constructed.
	As for the previous example, we observe a reduction of the theoretically stable timestep as given in \cref{tab:2DCEE} which we attribute again to the lack of \ac{ssp} guarantees.
	\begin{table}
		\def\arraystretch{1.2}
		\centering
		\begin{tabular}{c?{2}c|c|c|c|c}
			$S$ & $\Delta t$ & CFL & $N_t$  & $ \widetilde{\mathcal{M}} \cdot 10^{-15}$ & $\Delta t^3$ \\
			\Xhline{5\arrayrulewidth}
			30 & $5.92 \cdot 10^{-2}$ & $0.79$ & $338$ & $ 1.80 \cdot 10^{-12}$ & $2.08 \cdot 10^{-4}$ \\
			60 & $9.75 \cdot 10^{-2} $ & $0.65$ & $206$ & $3.21 \cdot 10^{-11}$ & $9.26 \cdot 10^{-4}$\\
			120 & $1.20 \cdot 10^{-1} $ & $0.40$ & $167$ & $6.50 \cdot 10^{-10}$ & $1.72 \cdot 10^{-3}$
		\end{tabular}
		\caption[Experimental data 2D compressible Euler equations]
		{Results for the simulation of the 2D compressible Euler equations \cref{eq:EulerEqs} with three many stage optimized methods constructed using \cref{alg:RootOpt}.
			The timestep $\Delta t$ is kept constant over the $N_t$ timesteps, with the possible exemption of the last timestep which may be reduced in order to match the desired final time $t_f = 20$.}
		\label{tab:2DCEE}
	\end{table}
	For the $L^\infty$ error in density $\Vert \rho^{(h)}(t_f, x, y) - \rho(t_f, x, y) \Vert_\infty$ second order convergence is observed, until spatial discretization errors start to become significant, cf. \cref{fig:2DCEE}.
	\section{Conclusions}
	\label{sec:Conclusions}
	In this work, a novel optimization approach for the generation of optimal stability polynomials for spectra of hyperbolic \ac{pdes} is devised.
	Parametrizing the stability polynomials in terms of the herein introduced pseudo-extrema offers both a numerically stable and readily interpretable representation.
	The optimization approach is motivated by examining the properties of the pseudo-extrema for the proven optimal stability polynomials of first and second order for disks.
	These findings are directly extended to strictly convex spectra for which an optimization procedure is developed in detail.
	For non-convex spectra, possible remedies are discussed and optimal stability polynomials for both convex and nonconvex spectra are presented.
	Stability polynomials with degrees larger than 100 are constructed for a range of classical hyperbolic \ac{pdes} that match the linear consistency requirements up to order three.
	
	Equipped with high degree stability polynomials in factorized form a possible choice for the actual Runge-Kutta methods is proposed.
	The construction of the numerical schemes centers around minimizing the propagation and amplification of round-off errors, which can become a concern for the many stage methods.
	Here, methods of second order are constructed which match their designed order of accuracy for both linear and nonlinear problems.
	For linear problems, only internal stability might limit the theoretically possible maximum timestep while for nonlinear problems the lack of the \ac{ssp} property is much more severe which spoils the effectiveness of the very high-stage methods.
	Nevertheless, the medium stage methods with up to, say, 64 stages are a promising candidate for stabilized methods as they allow an increase of the characteristic speeds in a simulation by a factor of four at a constant timestep compared to schemes constructed using a monomials-based approach.

	\section*{Data Availability}
	All data generated or analysed during this study are included in this published article and its supplementary information files.
	\section*{Acknowledgments}
	Funding by German Research Foundation (DFG) under grant number DFG-FOR5409.

	\section*{Declaration}
	The authors have no relevant financial or non-financial interests to disclose.

	\section*{Peer-Reviewed, Published Article}
	The peer-reviewed version of this article is available at \url{https://rdcu.be/dBAvS} with DOI \url{https://doi.org/10.1007/s10915-024-02478-5}.

	\appendix
	
	\section{Pseudo-Extrema of Shifted Chebyshev Polynomials}
	\label{sec:PECheby}
	Consider the shifted Chebyshev polynomials of first kind
	\begin{equation}
		\label{eq:ShiftedChebyshevPolyFirstKind}
		T_S\left(1 + z/S^2\right) = \cos \big( S \arccos \left( 1 + z/S^2 \right) \big) 
	\end{equation}
	providing the optimal first order stability polynomial for parabolic spectra \cite{franklin1958numerical, guillou1960domaine, dinsome1958}.
	It is well-known that the $S+1$ extrema of the shifted Chebyshev polynomials are on the $[-2S^2, 0]$ interval given by 
	\begin{equation}
		\label{eq:ChebyshevExtremePoints}
		x_j \coloneqq S^2 \left( \cos\left( \frac{j \pi}{S} \right) -1 \right), \quad j = 0, \dots, S
	\end{equation}
	where $T_S\left(1 + x_j/S^2\right)$ attains either $+1$ or $-1$.
	In particular, \cref{eq:ChebyshevExtremePoints} includes the $S-1$ critical points corresponding to $j = 1, \dots , S-1$ besides the extrema at the ends of the domain which correspond to $j = 0$ (right end) and $j = S$ (left end).
	
	We can rewrite the stability polynomial according to \cref{eq:DefinitionLowerDegree} as
	\begin{equation}
		\label{eq:LowerDegreeLHS}
		\widetilde{P}_{S-1}(z; \widetilde{\boldsymbol r}) = \frac{P_{S,1}(z; \widetilde{\boldsymbol r}) - 1}{z}
	\end{equation}
	and seek out to determine the pseudo-extrema $\widetilde{\boldsymbol r}$, i.e., the roots of $\widetilde{P}_{S-1}$.
	Apart from $z=0$ the roots of $\widetilde{P}_{S-1}$ are given by the roots of the nominator, i.e., 
	\begin{equation}
		\label{eq:EqChebyExtremePoints}
		P_{S,1}(z; \widetilde{\boldsymbol r}) - 1 = 0.
	\end{equation}
	For even degree $S$, $S/2+1$ roots of \cref{eq:EqChebyExtremePoints} are given by the Chebyshev extreme points where $T_S\left(1 + x_j/S^2\right)$ takes value $+1$:
	\begin{equation}
		\widetilde{r}_j = S^2 \left( \cos\left( \frac{ 2 j \pi}{S} \right) -1 \right), \quad j = 0, \dots, S/2.
	\end{equation}
	We recall that the nominator of \cref{eq:LowerDegreeLHS} $P_{S,1}(z) -1= T_S\left(1 + z/S^2\right) -1$ is a polynomial of degree $S$ in real coefficients and thus has, by the fundamental theorem of algebra, $S$ (possibly complex-conjugated) roots.
	As we have found already $S/2+1$ real roots it follows that the remaining $S/2 - 1$ roots have to be multiples of the already found ones.
	In particular, we have that the ``outer`` roots $\widetilde{r}_0 = -2S^2$ and $\widetilde{r}_{S/2} = 0$ have multiplicity one, and the ``interior`` roots $\widetilde{r}_j, j = 1, \dots S/2 -1$ have each multiplicity two.
	This follows from the fact that the $\widetilde{r}_j, j = 1, \dots S/2-1$ are a subset of the critical points of $P_{S, 1}(z)-1$, i.e., $\frac{\nid}{\nid z} \left(T_S\left(1 + z/S^2\right) -1 \right)\vert_{\widetilde{r}_j}= 0, \: j = 1, \dots S/2 -1$.
	
	Considering the roots of \cref{eq:LowerDegreeLHS} we have to be careful when it comes to $j = 0$, i.e., $x_j = 0$ since we divide by $\widetilde{r}_j = 0$.
	Recalling the introduction of the lower degree polynomial \eqref{eq:DefinitionLowerDegree} we have by construction $\widetilde{P}_{S-1}(0; \widetilde{\boldsymbol r}) = 1$.
	As a consequence, we have that the $S-1$ roots of the lower-degree polynomial $\widetilde{P}_{S-1}(z; \widetilde{\boldsymbol r})$ are given by 
	\begin{equation}
		\label{eq:PseudoExtremaParabolic}
		\widetilde{r}_j = S^2 \left( \cos\left( \frac{ 2 j \pi}{S} \right) -1 \right), \quad j = 1, \dots, S/2
	\end{equation}
	where $\widetilde{r}_{S/2}$ has multiplicity one and the remaining $S/2-1$ roots have multiplicity two.
	%
	%
	
	For odd degree $S$, we can follow a similar reasoning which yields that the pseudo-extrema in this case are given by 
	\begin{equation}
		\label{eq:PseudoExtremaOddLine}
		\widetilde{r}_j = S^2 \left( \cos\left( \frac{2j \pi}{S} \right) -1 \right), \quad j = 1, \dots, \frac{S-1}{2}
	\end{equation}
	which all have multiplicity two.
	This follows from the fact that for odd $S$, the extremum at the left end of the $\left[-2S^2, 0\right]$ interval of the Chebyshev polynomial is $-1$ while for even $S$, the left extremum takes value $+1$.
	
	While these findings are for the shifted Chebyshev polynomials admittedly trivial, they generalize neatly to the known first and second order accurate optimal stability polynomials of circular/disk-like spectra, as shown in \cref{sec:PEProvenOptimalStabPolys}.
	\section{Simulation Configurations with Strictly Convex Spectra}
	\label{app:SpectraSetup}
	\begin{enumerate}
		\item \label{item:Burgers} Burgers' Equation $\partial_t u + \partial_x \frac{1}{2} u^2 = 0$ is discretized using the \ac{dgsem} on $ \Omega = [0, 1]$ and Godunov flux \cite{OsherRiemann1984}. The periodic domain is discretized with $64$ elements on which the solution is reconstructed using third order local polynomials.
		As initial condition we supply $u(t_0, x) = 2 + \sin(2 \pi x)$ leading to a multi-valued solution at $t = \frac{1}{2 \pi}$.
		The spectrum is computed using the algorithmic differentiation capabilities of \texttt{Trixi.jl} and is displayed in \cref{fig:ConvexSpectraBurgers}.
		\item \label{item:ShallowWater} The Shallow Water equations with variable bottom topography $b(x)$ read in 1D 
		\begin{equation}
			\label{eq:ShallowWater1D}
			\partial_t \begin{pmatrix}
				h \\ hv 
			\end{pmatrix}
			+ \partial_x \begin{pmatrix} hv \\ hv^2 + \frac{g}{2}h^2 
			\end{pmatrix} + \begin{pmatrix}
				0 \\ gh \partial_x b
			\end{pmatrix}  = \boldsymbol 0.
		\end{equation}
		We discretize \cref{eq:ShallowWater1D} on $\Omega = [0, \sqrt{2}]$ using the \ac{dgsem} with flux differencing \cite{doi:10.1137/S003614290240069X, CHEN2017427}.
		In particular, for the first component of \cref{eq:ShallowWater1D} the surface flux is approximated using the Rusanov/Local Lax-Friedrichs flux and the second component via the flux proposed in \cite{FJORDHOLM20115587}.
		For the volume term we used the fluxes presented in \cite{WINTERMEYER2017200}.
		$\Omega$ is discretized with $32$ cells where we use third order local polynomials to reconstruct the solution.
		The initial condition is in this case a moderate discontinuity:
		\begin{equation}
			\begin{pmatrix}
				h(t_0, x) \\ v(t_0, x)
			\end{pmatrix} = 
			\renewcommand*{\arraystretch}{2.5}
			\begin{pmatrix} - b(x) + 
				\begin{cases} 3.25 & \vert x- 0.7 \vert > 0.5 \\ 4 & \text{else}	\end{cases}	\\ \qquad \qquad \: \: \: \begin{cases} 0 & \vert x- 0.7 \vert > 0.5 \\ 0.1882 & \text{else}	\end{cases}	\end{pmatrix}
		\end{equation}
		with sinusoidal bottom topography $b(x) = \sin(x)$.
		The corresponding spectrum is displayed in \cref{fig:ConvexSpectraShallowWater}.
		\item \label{item:MHD} The ideal compressible \ac{mhd} equations considered here are of form
		\begin{equation}
			\partial_t \begin{pmatrix}
				\rho \\ \rho v_1 \\ \rho v_2 \\ \rho v_3 \\ \rho e \\ B_1 \\ B_2 \\ B_3 
			\end{pmatrix}
			+ \partial_x \begin{pmatrix}
				\rho v_1 \\ 
				\rho v_1^2 + p + e_B - B_1^2 \\
				\rho v_1 v_2 - B_1 B_2 \\ 
				\rho v_1 v_3 - B_1 B_3 \\
				v_1(e_v + \gamma (\rho e - e_v - e_B) + 2 e_B) - B_1(v_1 B_1 + v_2 B_2 + v_3 B_3) \\
				0 \\
				v_1 B_2 - v_2 B_1 \\
				v_1 B_3 - v_3 B_1 
			\end{pmatrix}
			= \boldsymbol 0
		\end{equation}
		with $e_B \coloneqq 0.5 (B_1^2 + B_2^2 + B_3^2), e_v \coloneqq \rho (v_1^2 + v_2^2 + v_3^2), p = (\gamma - 1) (\rho e - e_v - e_B)$.
		The solution is represented with $p=1$ local polynomials and we employ the Rusanov/Local Lax-Friedrichs flux \cite{RUSANOV1962304}.
		The periodic domain $\Omega = [0, 1]$ is again discretized using $32$ cells.
		For the initial condition we choose the Alfv\'{e}n wave \cite{TOTH2000605}
		\begin{equation}
			\begin{pmatrix}
				\rho (t_0, x) \\ v_1 (t_0, x) \\ v_2 (t_0, x) \\ v_3 (t_0, x) \\ p (t_0, x) \\ B_1 (t_0, x) \\ B_2 (t_0, x) \\ B_3 (t_0, x)
			\end{pmatrix} = \begin{pmatrix}
				1 \\ 0 \\ 0.1 \sin (2 \pi x) \\ 0.1 \cos(2 \pi x) \\ 0.1 \\ 1 \\ 0.1 \sin(2 \pi x) \\ 0.1 \cos(2 \pi x)
			\end{pmatrix}
		\end{equation}
		and set $\gamma = \frac{5}{3}$.
		The spectrum for this particular configuration is shown in \cref{fig:ConvexSpectraMHD}.
		\item \label{item:2D_CEE} The compressible Euler equations in two spatial dimensions are given by 
		\begin{equation}
			\partial_t
			\begin{pmatrix}
				\rho \\ \rho v_1 \\ \rho v_2 \\ \rho e
			\end{pmatrix}
			+
			\partial_x
			\begin{pmatrix}
				\rho v_1 \\ \rho v_1^2 + p \\ \rho v_1 v_2 \\ (\rho e +p) v_1
			\end{pmatrix}
			+
			\partial_y
			\begin{pmatrix}
				\rho v_2 \\ \rho v_1 v_2 \\ \rho v_2^2 + p \\ (\rho e +p) v_2
			\end{pmatrix}
			= \boldsymbol 0
		\end{equation}
		which are represented on $\Omega = [-1, 1]^2$ with $8$ cells per direction using second order local polynomials.
		Here, $p = (\gamma - 1) (\rho e - 0.5 \rho (v_1^2 + v_2^2))$ and $\gamma = 1.4$.
		For the flux we use the \ac{hllc} flux \cite{toro1994restoration} and supply again periodic boundaries.
		The fields are initialized as 
		\begin{equation}
			\begin{pmatrix}
				\rho(t_0, x) \\ v_1 (t_0, x) \\ v_2 (t_0, x) \\ p (t_0, x)
			\end{pmatrix}	
			= \begin{pmatrix}
				2 + 0.1 \sin \big(\pi (x + y)\big) \\ 1 \\ 1 \\ \big[2 + 0.1 \sin \big(\pi (x + y)\big)\big]^2
			\end{pmatrix}.
		\end{equation}
		The spectrum corresponding to this initial state and semidiscretization is displayed in \cref{fig:ConvexSpectra2D_CEE}.
	\end{enumerate}
	%
	

	\bibliographystyle{elsarticle-num-names} 
	\bibliography{references.bib}
\end{document}